\documentclass{IET}

\usepackage{multirow}
\usepackage[table,xcdraw]{xcolor}
\usepackage[color=yellow]{todonotes}

\usepackage{arydshln}
\usepackage{subfigure}
\usepackage{graphicx}
\usepackage{caption}

\usepackage[noadjust]{cite}

\usepackage{algorithm}
\usepackage{algpseudocode}
\usepackage{fancyhdr}

\fancypagestyle{firststyle}
{
   \fancyhf[C]{This paper is a postprint of a paper submitted to and accepted for publication in  IET Generation, Transmission \& Distribution and is subject to Institution of Engineering and Technology Copyright. The copy of record is available at IET Digital Library}
   \fancyfoot[C]{ \thepage\ }
}
\begin{document}

\title{A hybrid approach for planning and operating active distribution grids}

\author{Stavros Karagiannopoulos}
\author{Petros Aristidou}
\author{Gabriela Hug}
\affil{EEH - Power Systems Laboratory, ETH Zurich, Physikstrasse 3, 8092 Zurich, Switzerland \\
E-mail: \{karagiannopoulos, aristidou, hug\}@eeh.ee.ethz.ch}

\abstract{This paper investigates the planning and operational processes of modern distribution networks (DNs) hosting Distributed Energy Resources (DERs). While in the past the two aspects have been distinct, a methodology is proposed in this paper to co-optimize the two phases by considering the operational flexibility offered by DERs already in the planning phase. By employing AC Optimal Power Flow (OPF) to analyse the worst-case forecasts for the load and distributed generator (DG) injection, the optimal set-points for the DERs are determined such that the network's security is ensured. From these results, the optimized individual characteristic curves are then extracted for each DER which are used in the operational phase for the local control of the devices. The optimized controls use only local measurements to address system-wide issues and emulate the OPF solution without any communication. Finally, the proposed methodology is tested on the Cigre LV benchmark grid confirming that it is successful in mitigating with acceptable violations over- and under-voltage problems, as well as congestion issues. Its performance is compared against the OPF-based approach and currently employed local control schemes.}
\maketitle
\thispagestyle{firststyle} 

\section{Introduction}\label{sec:intro}

Some of the most noticeable developments foreseen in power systems involve active Distribution Networks (DNs). In the future, DNs are expected to host a large percentage of the Renewable Energy Sources (RES) and, along with other small generators connected at the distribution level, to supply a growing share of the total demand \cite{Lopes2007,EuropeanCommission2015}. These units, in combination with other Distributed Energy Resources (DERs) such as electric vehicles, energy storage systems and flexible loads, will amplify the role of distribution grids, allowing them to provide power and ancillary services to higher voltage levels. However, such new developments bring significant challenges in both the DN planning and operation stages. 

In the operation stage, the increased penetration of RES can lead to over-voltages, reverse flows, and thermal overloads at times of high DG generation and low consumption. In response to these problems, several DN control schemes have been proposed for the dispatch of DERs~\cite{Joos2000}. The available  control schemes differ with respect to two aspects: First, depending on the type and level of communication infrastructure required, the DN control schemes can be categorized as local (without communication), distributed (limited, peer-to-peer communication), and centralized (full communication). Another important distinction among these methods concerns the off-line knowledge about the system topology, parameters, or operation (historical or forecasted data), to the controllers. An overview of such schemes is presented in Section~\ref{sec:overview}.

In the planning stage, the RES generation and the load growth and consumption profiles need to be accurately forecasted. The complexity of planning methods increases even more when DERs become active elements, i.e. if their active and reactive power injections become a function of the DN operating conditions. It is very important to include this active behaviour of DERs in the planning phase to ensure compatibility between the two stages (planning and operation). 

Current procedures employed by Distribution System Operators (DSOs) clearly separate the planning stage from the operation stage~\cite{Cigre2014,SmartPlanningWP2.5}. Based on critical snapshots of the forecasted RES generation and load consumption, a series of off-line load flow calculations is employed to ensure that the worst case conditions in terms of voltage thresholds, component loading, and other security constraints, can be satisfied. Using the outcome of this analysis, DSOs make decisions on reinforcements in their system to alleviate any violations. The flexibility of DERs is usually not considered, unless the existing grid codes dictate an active DER response. Current regulations concern local reactive power schemes, based on the unit's active power or terminal voltage.

Such local approaches have the advantage of being simple and straightforward to implement. However, the units are programmed to react the same, irrespective of the system topology or their exact location inside the network. In the near future, centralized solutions will become available due to the installation of smart metering devices and dispatchable DERs, and might lead to more efficient DN operation.

For this reason, optimal planning of active distribution systems revisits traditional planning methods while considering the new control  and communication capabilities of modern DNs. This area has seen a lot of development during the last decade \cite{Paiva2005,Mohtashami2016,Naderi2012,Gan2011,Pilo2013,Georgilakis2015}. Among the several approaches proposed, the dominant ones are optimization-based (MILP~\cite{Paiva2005},  MINLP~\cite{Mohtashami2016}, etc.), and heuristic (genetic and evolutionary algorithms~\cite{Naderi2012}, etc.). Other methods are based on statistical analysis to obtain the optimal LV network design, while considering the life cycle cost of equipment, load density, substation requirement, and grid type inputs~\cite{Gan2011}, or probabilistic approaches combined with Monte Carlo simulations to capture uncertainties and risks~\cite{Pilo2013}. A more detailed review on models, methods, and future of distribution system planning can be found in \cite{Georgilakis2015}.

The transition from the currently prevailing local approach to a fully central operation, as well as the necessity for higher interdependence of the planning and operation stages have not been studied adequately. In this paper, we present a methodology for planning and operating modern DNs with active elements which postpones investments in communication infrastructure and new equipment. The method co-optimizes these two stages by transferring information from the planning phase to the operation phase. First, an offline centralized approach is used where the optimal behaviour of the controlled DERs is obtained based on load consumption and RES generation forecasts from the planning phase. Then, from these results, individual local characteristics are derived for each unit to be used in the real-time operation. In this way, each DER has a unique response to local measurements (active power generation and terminal voltage), based on its location in the DN and available offline information.

The remainder of the paper is organized as follows. First, Section~\ref{sec:overview} provides an overview over present planning and operating aspects of active distribution grids. Then, the proposed methodology is analysed in Section~\ref{sec:Method}, while Section~\ref{sec:CaseStudies} introduces the considered case studies and presents the simulation results. Finally, conclusions are outlined in Section~\ref{sec:Conclusions}.

\section{Overview of distribution grid planning and operation stages}\label{sec:overview}

In this section, an overview of current practices and trends in DN planning and operation are summarized.

\subsection{Planning stage}

Recent surveys~\cite{Cigre2014,SmartPlanningWP2.5}, show that the introduction of DERs and new load types in DNs is not  appropriately reflected in present planning procedures. On the contrary, planning techniques are still based on a conventional approach~\cite{SmartPlanningWP2.5}.  
First, all relevant data, including historical weather and electricity demand, are gathered. 
Next, each DSO develops forecasts of the future demand and generation, based on which the planning decisions {\color{black} will be made}. Then, the DSOs perform deterministic load flow analyses of the worst-cases scenarios, leading to - often conservative - decisions regarding grid reinforcements. Different feasible solutions are identified according to the DSO's targets and compliance with policies and international standards/regulations. Finally, in the last stage, the available options are evaluated based on different criteria, such as  lifetime network costs, implementation complexity, or reliability targets.

These conventional planning techniques do not include the capabilities of the new technologies. They are clearly separated from the operation stage and thus, new planning guidelines are needed to fully incorporate the new active technologies.  In the future, data gathered in the planning stage should also be used to steer the design of controls in the real-time operation, something that does not occur today.

\subsection{Operation stage}

Triggered by financial incentives -at least partially- due to subsidies as well as the desire to render the energy supply system more sustainable, RES, energy storage systems, and small CHP plants, comprise already a significant share in the generation mix of most countries~\cite{EuropeanCommission2015}. In the past, these units were operating with a power factor  close to one and were regarded as negative loads injecting pure active power. However, with the recent introduction of new grid-codes and regulations, DERs are called upon to provide ancillary services to the grid, particularly for voltage support.

Based on the communication infrastructure available in the DN, the control schemes can be classified as local, distributed, or centralized.

\subsubsection{Local control schemes} 

These schemes are widely analysed in the literature~\cite{Weckx2014,Demirok2011,Demirok2010EvalVOLT,Ellis2012} and have been already incorporated into grid-codes and standards~\cite{VDE} to mitigate voltage problems. Their main feature is that the DERs react to purely local measurements according to some predefined rules. No communication or remote measurements are available for these controllers.

Some of the most common rules are for the DGs to provide reactive power compensation based on their terminal voltage measurements, i.e. $Q=f(V)$~\cite{Kotsampopoulos2013}, or their voltage and active power injection measurements, i.e. $Q=f(V,P)$~\cite{Weckx2014}. Another approach is to adjust the power factor according to their active power injection, i.e $cos(\phi)= f(P)$~\cite{ItalianNorm2014,VDE}. When reactive power compensation is not sufficient to mitigate the voltage issues, active power curtailment can also be employed locally based on voltage measurements, i.e. $P_{curt}=f(V)$~\cite{Weckx2014,Tonkoski2011}. Other more simplified approaches include a fixed power factor mode or fixed reactive power absorption.  

\begin{figure}
  \centering
      \subfigure[]{\includegraphics[width=0.7\textwidth]{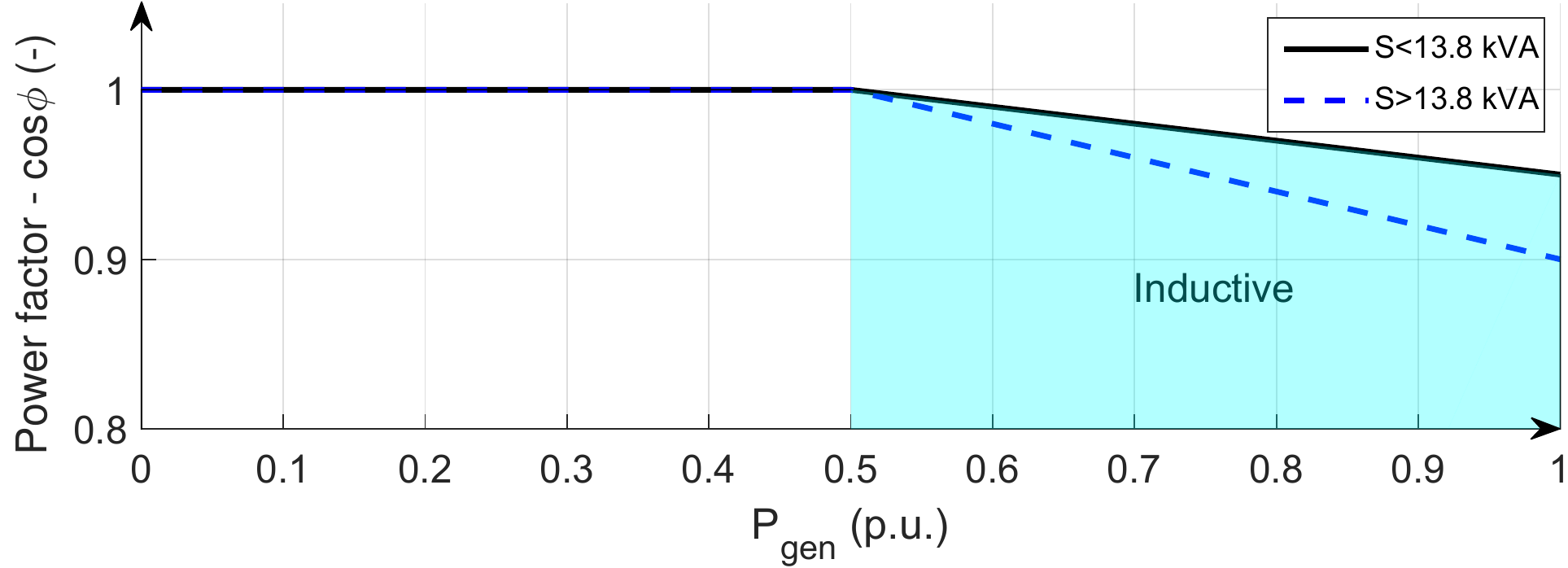}\label{pfchar}}
      \subfigure[]{\includegraphics[width=0.7\textwidth]{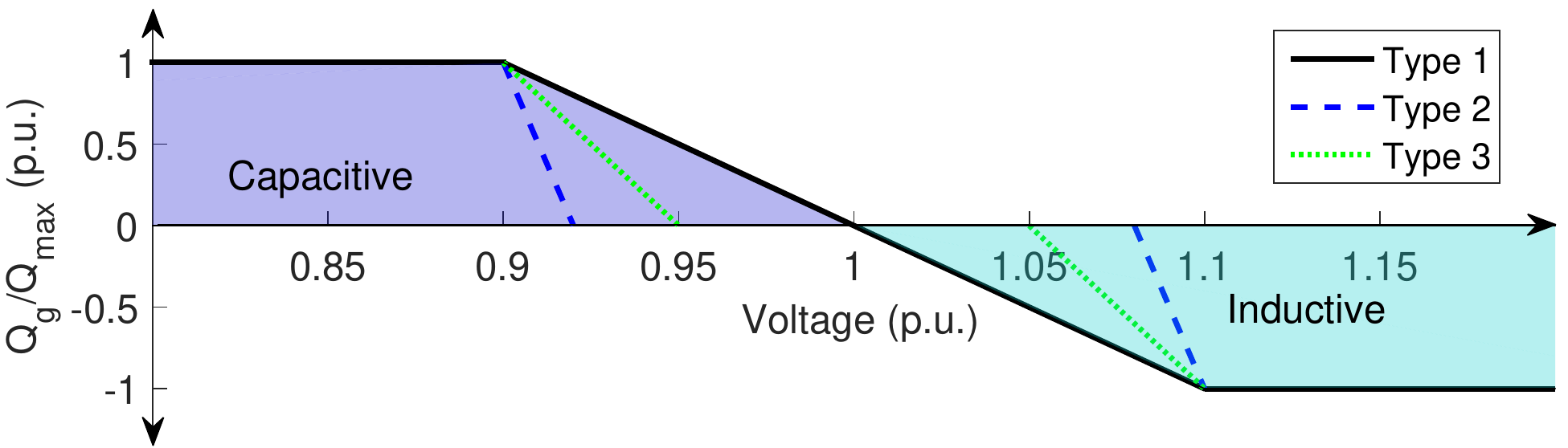}\label{localschemes}}
      \subfigure[]{\includegraphics[width=0.7\textwidth]{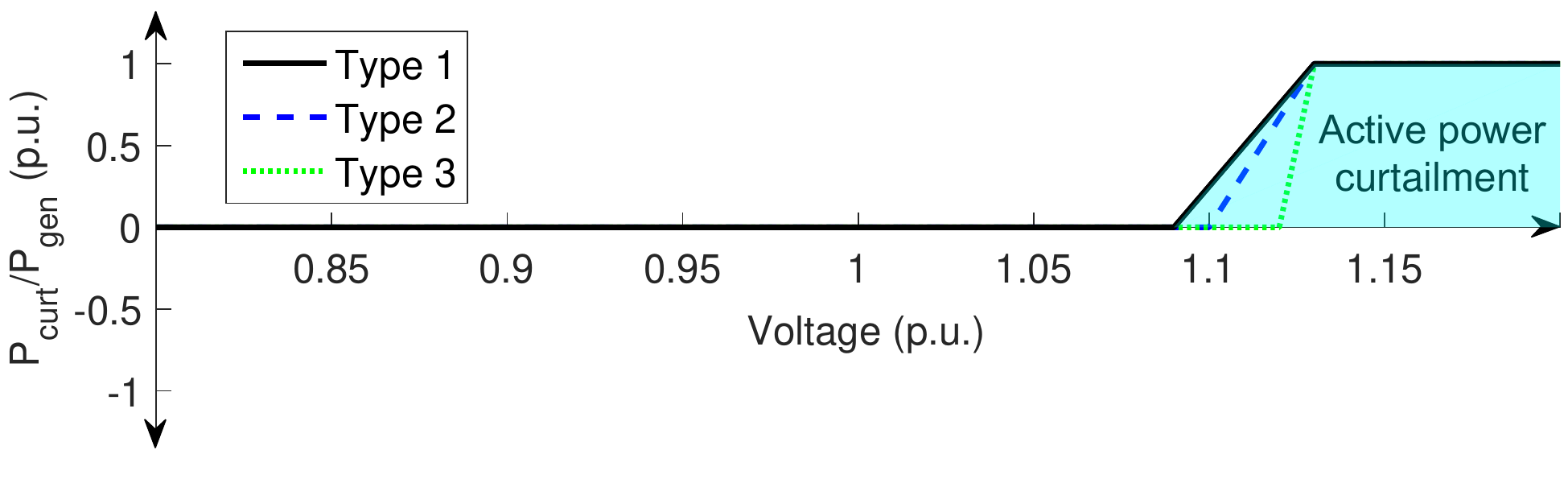}\label{localschemes_P}}
      \subfigure[]{\includegraphics[width=0.65\textwidth]{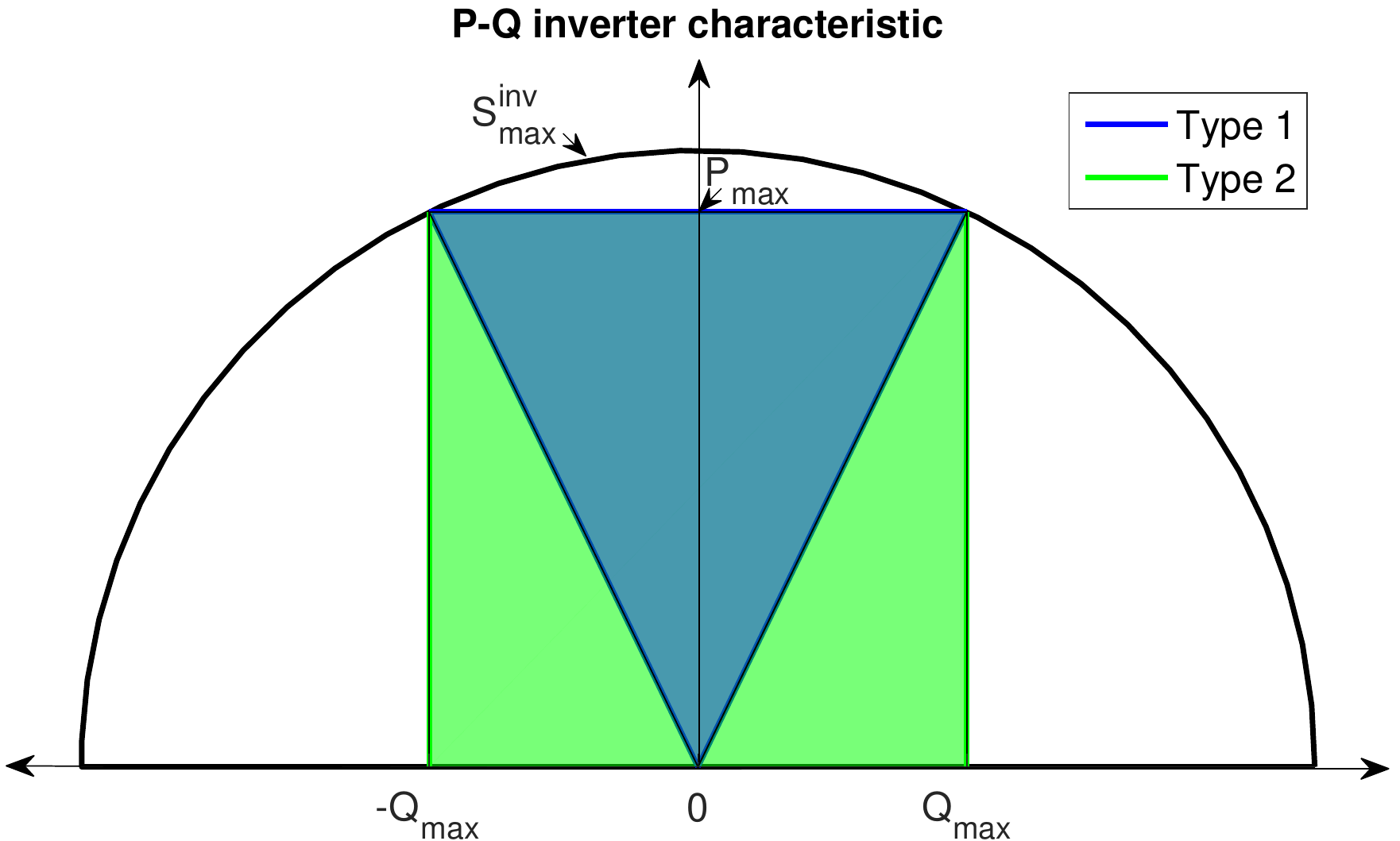}\label{PQcap}}
      \caption{  a) Default power factor inverter characteristic in Germany~\cite{VDE}.\\ b) Alternative local reactive power control schemes.\\
      c) Alternative local active power curtailment schemes.\\ d) Steady state $P-Q$ inverter capability curve.}
  \end{figure}

Figure~\ref{pfchar} shows the present $cos\phi=f(P)$ characteristic curve implemented in Germany. As soon as the DGs inject more than half of their maximum power, they have to adjust their power factor. Figures~\ref{localschemes} and~\ref{localschemes_P} show different  versions of local reactive and active power control, namely $Q=f(V)$ and $P_{curt}=f(V)$. Low local voltages impose a capacitive behaviour and production of reactive power to increase the voltage. The opposite behaviour is observed at high voltages. The different types define the level of responsiveness to voltages at different operating points. E.g., Type 1 results in a continuous consumption/injection of reactive power, while Types 2 and 3 allow for a dead-band, in which the inverter control is inactive. Regarding active power curtailment, Fig.~\ref{localschemes_P} shows similar types of responsiveness to voltages. Type 1 initiates the curtailment before the voltage reaches the threshold of $1.1$ p.u., whereas the other two start curtailing after that. 

\subsubsection{Distributed control schemes}

These schemes employ communication between neighbouring DERs in order to achieve a better operation than the pure local laws. For example, reactive compensation can be provided based not only on the unit's terminal voltage measurements but also of neighbouring units, i.e. $Q=f(V_1,\ldots,V_n)$ where $n$ is the number of neighbouring DGs communicating their voltage measurements~\cite{Olivier2015}. The extra information available to the DGs allows to employ more sophisticated schemes and coordinate their resources to solve problems not possible with purely local schemes.

\subsubsection{Centralized control schemes} 

Finally, there are schemes which work in a centralized manner. On the one hand, rule-based centralized methods focus on controlling devices based on real-time measurements and some predefined rules derived from the historic knowledge of the system~\cite{Roytelman2000CoordinatedSystemsb},~\cite{WilliamsDistributionFlow},~\cite{MatsLarsson2000CoordinatedSystems}. On the other hand, optimization-based approaches rely on full communication to aggregate the measurements to a central entity where network-level optimization techniques are employed to compute the system-wide optimal settings for the controlled DERs. These methods are gaining great importance lately \cite{KaragiannopoulosGM,Fortenbacher2016,DallAnese2015,DallAnese2014} due to advances in information technology and new theoretical developments in approximations of the nonlinear AC power flow equations~\cite{Fortenbacher2016,Low2014} allowing faster computation times. For example, in~\cite{DallAnese2014} the optimal active and reactive power set-points of DGs are computed according to various objectives (such as minimization of losses). Similarly, in~\cite{KaragiannopoulosGM}, an OPF formulation is described by which the reactive and active power of the PV inverters are controlled centrally to mitigate over-voltages and to defer investments.

\subsubsection{ICT requirements, grid controllability, and optimal operation} 
Most of the newer-type inverters are able to get local measurements, change their power factor, and limit the active power output, through some primary controller. These functionalities (i.e., the ability to get local measurements and change the P,Q set-points) are necessary for any type of control scheme - either local or centralized - that relies on controlling the DERs.

Centralized control schemes offer increased observability and controllability of the grid, thus allowing us to operate the DN close to optimal. However, they require a reliable two-way communication infrastructure between the measurement devices, the dispatchable DERs, and a central unit, as well as knowledge of the grid topology and parameters. Unfortunately, such infrastructure does not exist yet in most DNs due to the costs involved. Moreover, centralized schemes are often more complex than the other methods and rely on the solution of sometimes computationally intensive optimization problems.

For these reasons, the most prevailing approaches today are based on local control, where each DG is designed to respond only to local measurements of voltage or power. These methods are cheap and simple, and the lack of a central control unit and communication increases the robustness and reliability of the scheme. Nevertheless, based only on local measurements and without any real-time coordination, these schemes are sub-optimal compared to a centralized scheme.

Finally, when limited communication is available between some of the units, distributed schemes can be used to achieve a more optimal operation than the local schemes. Still, if the communication and system information is not complete, these methods are not as optimal as in the centralized case. However, they are more robust and come at a lower cost and complexity.

\subsubsection{Inverter-based DER capabilities}\label{inverters}
Since the majority of DERs installed today in low-voltage DNs are inverter-based, an important part of designing any control scheme is to understand their capabilities. Reference~\cite{Kotsampopoulos2013} provides a review over recent requirements for voltage support in Europe and USA. Figure~\ref{PQcap} shows the steady-state $P-Q$ capability curve of one such inverter under constant voltage, where $S_{max}^{inv}$ is the inverter nominal apparent power.

Based on current practices, the inverter is usually over-sized compared to the installed capacity of the DER ($P_{max}$), in order to offer reactive power capability even when the DER is injecting its maximum active power. When only active power curtailment is allowed, the operating point can move along the $P$ axis up to $P_{max}$. Regarding reactive power capability, two variants are shown: First, the most common case defined in standards such as~\cite{VDE}, is based on a minimum power factor requirement, forming a triangle representing the feasible area. 
The second variant (Type~2) can be found in recent standards, e.g.~\cite{ItalianNorm2014}. Here, the maximum reactive power is independent of the active power injection. 

\section{Co-optimization of planning and operating stages}\label{sec:Method}

In this section, the proposed method is explained in detail. The main idea is to use information from the planning stage in combination with an offline centralized control, in order to design a real-time local control scheme that resembles the optimal centralized control scheme without the need of any communication links. In this approach, system-wide information about the system topology and optimal behaviour is incorporated into the local controllers, thus customizing their characteristics to the specific DN requirements without the need to invest in additional communication equipment. Thus, a co-optimizing framework between the planning and operation stages is defined. 

\begin{figure}
	\centering
	\includegraphics[width=1\textwidth]{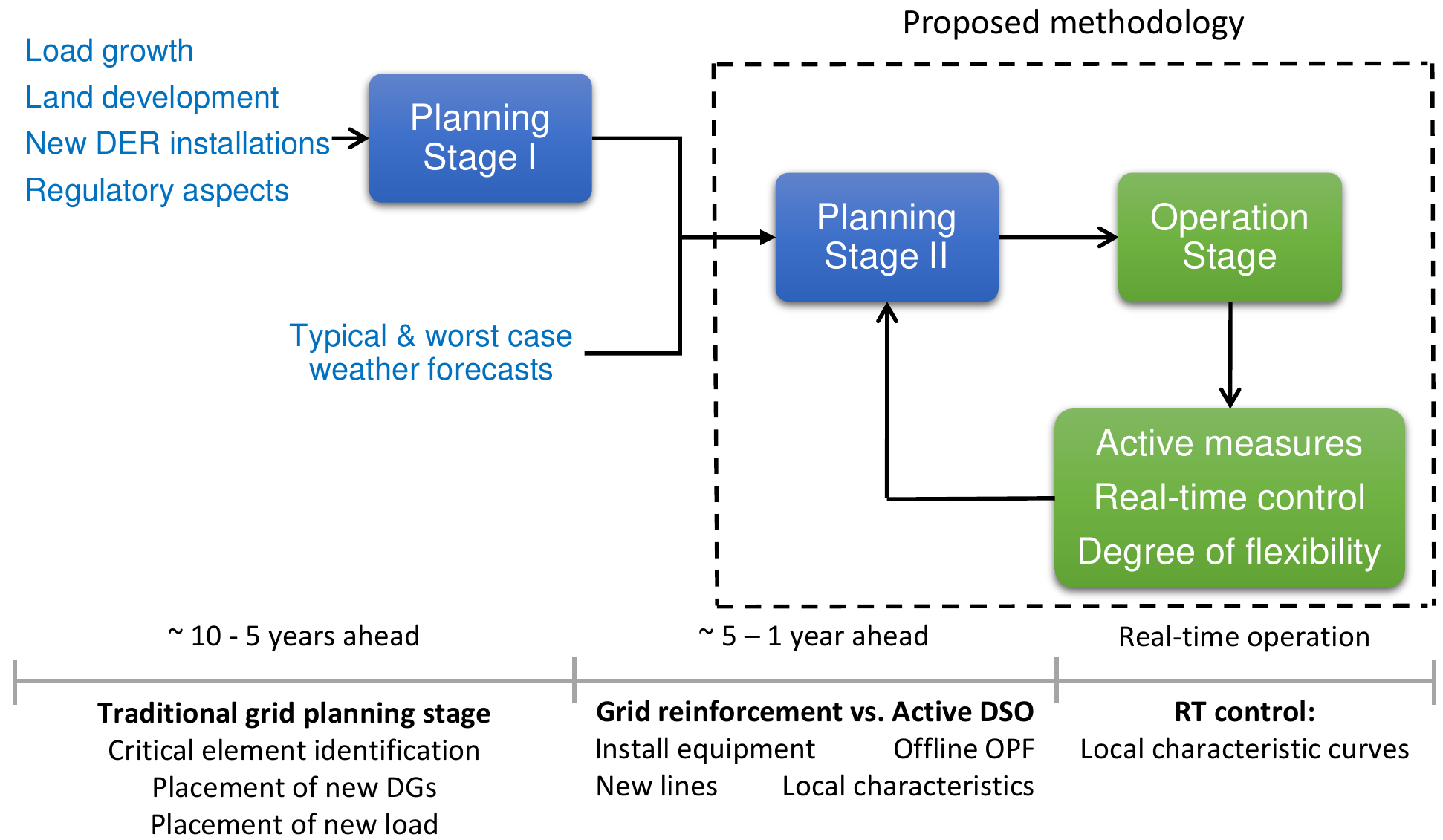}
	\caption{Proposed joint planning and operation stages.}
	\label{fig4}
\end{figure}

Figure~\ref{fig4} outlines the proposed framework. Planning Stage I resembles the traditional planning process, where the DSO decides on the general expansion needs based on the location of new DER installations, regulatory aspects and development plans.

In the Planning Stage II, the DSO takes the final decision on whether to invest in new equipment, e.g. lines/transformers/etc., or whether it will have adequate flexibility (active measures) in the operation stage postponing the investments. Here, the DSO has to solve a planning problem, i.e. should it build a new line or not, considering also operational aspects. An economic evaluation similar to~\cite{KaragiannopoulosGM} will assess the trade-offs between installing hardware with burdens on capital expenditure and using control measures increasing the operational costs. The outcome of this step focuses on a set of proposed control approaches for the flexible units. Thus, there is a stronger interdependence between the stages of planning and operation compared to the past. In case the economic evaluation indicates that the investment in hardware is more meaningful, the DSO can proceed installing new equipment to avoid constraint violations. The proposed control approaches are considered to be a valid solution to the planning problem if the active units can relieve the grid of its future challenges, by using their flexibility and controllability. By understanding the benefits and opportunities of the operation stage, the DSOs could make use of such an approach in order to take cost efficient decisions. For example, costly investments can be postponed, provided that adequate flexibility is offered in the operation stage by controlling DGs. 

\subsection{Offline optimal power flow}

By running an offline OPF based on forecasts, we derive the optimal behaviour of the units as if we had full grid observability and controllability. This problem is formulated as follows:

\paragraph{Objective function}
The objective of the DSO is to minimize the losses of the grid and the use of active measures which come with some cost. Thus, the objective function reads
\begin{align}
	\min_{u}  \,\,  (\textrm{losses} + c^{T} \cdot u) &\label{eq:theta_slack}
\end{align}
where vector $c$ represents the operational costs associated with the activated control measures vector $u$. Modern DSOs may control the active/reactive output of DERs, activate some flexible loads, or modify the setpoints of storage units.
	
In more detail, the objective function in presence of active power curtailment and reactive power control is
	
\begin{align}
	\min_{P_{l},P_{curt},Q_{ctrl}}  \,\, \sum \limits_{t=1}^{N_{hor}} (c^{T}_{el} \cdot P_{l}(t) + c^{T}_{P} \cdot P_{curt}(t) + c^{T}_{Q} \cdot Q_{ctrl}(t))  &\label{eq:theta_slack2}
\end{align}
where at time interval $t$:
\begin{itemize}
\item $P_{l}(t)=\sum \limits_{k=1}^{n_{B}} \sum \limits_{m=k+1}^{n_{B}} (|P_{km}(t) + P_{mk}(t)|)$ is the sum of the losses on the DN lines and transformers;
\item $P_{curt}(t) = \sum \limits_{i=1}^{n_{DG}} (P_{max}(i,t) - P_{g}^{\textrm{f}}(i,t))$ is the sum of the curtailed power over all DG units;
\item $P_{max}(i,t)$ is the maximum active power the DG $i$ can inject and $P_{g}^{\textrm{f}}(i,t)$ is its actual forecasted infeed;
\item $Q_{ctrl}(i,t)$ is the reactive power set-point for DG $i$;
\item $c^{T}_{P}$ represents the cost of curtailing active power and $c^{T}_{Q}$ the cost of using reactive power. Setting $c^{T}_{Q} \ll c^{T}_{P}$ allows prioritizing the use of reactive power control.
\end{itemize}

To optimize the use of reactive power, we minimize the overall use of reactive power. The corresponding objective function is defined as 
\begin{align}
	Q_{ctrl}^{1}(t)=\sum \limits_{i=1}^{n_{DG}} | Q_{g}^{\textrm{f}}(i,t)| 
    \label{Q1}
\end{align}     
With this formulation, the DSO tries to minimize the reactive power $|Q_{g}^{\textrm{f}}(i,t)|$ in the planning stage, while keeping the power factor of DGs close to one. By doing so, the DSO maintains more flexibility for the operation stage, where potential real-time control may be needed.

\paragraph{Constraints}
	The OPF formulation includes the power balance equations at every node and time step $t$ as given by
\begin{align}
	P_{\textrm{inj}}^{\textrm{f}}&=P_{\textrm{g}}^{\textrm{f}} - P_{\textrm{l}}^{\textrm{f}} \label{eq:node_balance_P}\\
    Q_{\textrm{inj}}^{\textrm{f}}&=Q_{\textrm{g}}^{\textrm{f}} - Q_{\textrm{l}}^{\textrm{f}}  \label{eq:node_balance_Q}   
\end{align} 
where $(P_{\textrm{g}}^{\textrm{f}},Q_{\textrm{g}}^{\textrm{f}})$ are the active and reactive power injection forecasts of the DGs, $(P_{\textrm{l}}^{\textrm{f}},Q_{\textrm{l}}^{\textrm{f}})$ the active and reactive forecasted node demand and  $(P_{\textrm{inj}}^{\textrm{f}},Q_{\textrm{inj}}^{\textrm{f}})$ the corresponding net node injections into the distribution lines. 

The voltage constraints at every node  are given by  
\begin{align}
	& V_\textrm{min} \leq V_\textrm{i} \leq V_{\textrm{max}}  \label{eq:voltage_lim}
\end{align}
where $V_{\textrm{min}}$ and $V_{\textrm{max}}$ are the upper and lower acceptable voltage limits. 
Similarly, the thermal limits of the distribution lines are imposed by 
\begin{align}
	& 0\leq |S_{\textrm{i,j}}^{\textrm{f}}| \leq {S_{\textrm{i,j}}^{\textrm{max}}} \label{eq:thermal_lim2}
\end{align}
where $S_{\textrm{i,j}}^{\textrm{f}}$ is the forecasted apparent power flowing through the line connecting nodes $i$ and $j$, and $S_{\textrm{i,j}}^{\textrm{max}}$ the value corresponding to its upper thermal limit.

Furthermore, the DG limits are given by
\begin{align}
	P_{\textrm{g}}^{\textrm{min}} &\leq P_{\textrm{g}}^{\textrm{f}} \leq P_{\textrm{g}}^{\textrm{max}} \label{eq:PV_prod2}\\
	Q_{\textrm{g}}^{\textrm{min}} &\leq Q_{\textrm{g}}^{\textrm{f}} \leq Q_{\textrm{g}}^{\textrm{max}}
\end{align}
where $P_{\textrm{g}}^{\textrm{min}}$, $Q_{\textrm{g}}^{\textrm{min}}$, $P_{\textrm{g}}^{\textrm{max}}$ and $Q_{\textrm{g}}^{\textrm{max}}$ are the upper and lower limits for active and reactive generation. These limits vary depending on the type of the DER and the control schemes implemented. 
	
For {\color{black}RESs}, $P_{\textrm{g}}^{\textrm{max}}$ is limited by the maximum power they can produce. Moreover, it is usual for small inverter-based generators to have limitations on the power factor they can operate at. These limits are either technical or defined by the grid codes. For example, in case of having a feasible area of Type 1 described in section \ref{inverters}, the reactive power limit of \eqref{eq:PV_prod2} is modified to
\begin{align}
	 -\textrm{tan}(\phi_{\textrm{max}}) \cdot P_{\textrm{g}}^{\textrm{f}} \leq  Q_{\textrm{g}}^{\textrm{f}} \leq \textrm{tan}(\phi_{\textrm{max}}) \cdot P_{\textrm{g}}^{\textrm{f}} \; &\label{eq:PV_prod2Q}
\end{align}
where $\textrm{cos}\phi_{max}$ is the maximum allowed power factor.

Regarding Type 2 (see \ref{inverters}), we have a larger capability area for the reactive power. In this case the limit of \eqref{eq:PV_prod2} is modified to
\begin{align}
 	-\textrm{tan}(\phi_{\textrm{max}}) \cdot P_{\textrm{g}}^{\textrm{max}} \leq  Q_{\textrm{g}}^{\textrm{f}} \leq \textrm{tan}(\phi_{\textrm{max}}) \cdot P_{\textrm{g}}^{\textrm{max}} \; &\label{eq:PV_prod2Qb}
 \end{align}

Finally, in the case of RES with fixed unit power factor, the constraints become 
\begin{align}
   P_{\textrm{g}}^{\textrm{f}} &= P_{\textrm{g}}^{\textrm{max}} \label{eq:PV_prod2a_P}\\
   Q_{\textrm{g}}^{\textrm{f}} &= 0 \label{eq:PV_prod2a_Q}
\end{align}

\subsection{Derivation of inverter characteristics}\label{sec:deriv_char}

In this section, we present the procedure of deriving the individual characteristic curves for each DER. These curves dictate the real-time response of the units to local measurements and can be programmed on the DERs using lookup tables. By including the information from the offline optimization to customize the local controllers, the topology of the system as well as the forecasted data are implicitly taken into account in the real-time control.

Extracting the characteristic curves from the OPF-based solution depends on the type of local controller to be used by the DERs. For instance, a $Q=f(V)$ curve can be extracted considering the voltages as well as the reactive power injection/consumption of the OPF-based solution. 

The authors of~\cite{Demirok2010EvalVOLT} propose parametrized $Q=f(V)$ curves, adjusting mainly the deadband of the curves based on the impedance that each inverter sees from the grid. The minimum and maximum values of $V=0.9$ p.u and $V=1.1$ p.u., respectively, are regarded as fixed and the main focus is to have different deadband width around the fixed value of $V=1$ p.u.. In this paper, we use an off-line OPF approach which allows us to handle better all kinds of constraints and various types of objectives. Furthermore, the shape of the proposed characteristic curves can be composed of multiple levels, providing more flexibility compared to the ones currently implemented.

The authors of~\cite{Weckx2014} use a centralized optimization algorithm which frequently optimizes parameters of local controllers using real-time measurements. The focus is restricted to the unbalanced operational stage, using the open-loop control scheme $Q=f(P)$ and $P_{curt}=f(P)$. Furthermore, they use a linearized voltage model, in contrast to the more accurate non-linear AC power flow equations. In this paper, we follow a closed-loop $Q=f(V,P)$ scheme, showing better responsiveness to voltages. Furthermore, our OPF calculations are performed off-line using forecasts and worst-case profiles, assuming that there is no centralized control in real-time. 

\begin{figure}
  \centering
      \subfigure[]{\includegraphics[width=0.8\textwidth]{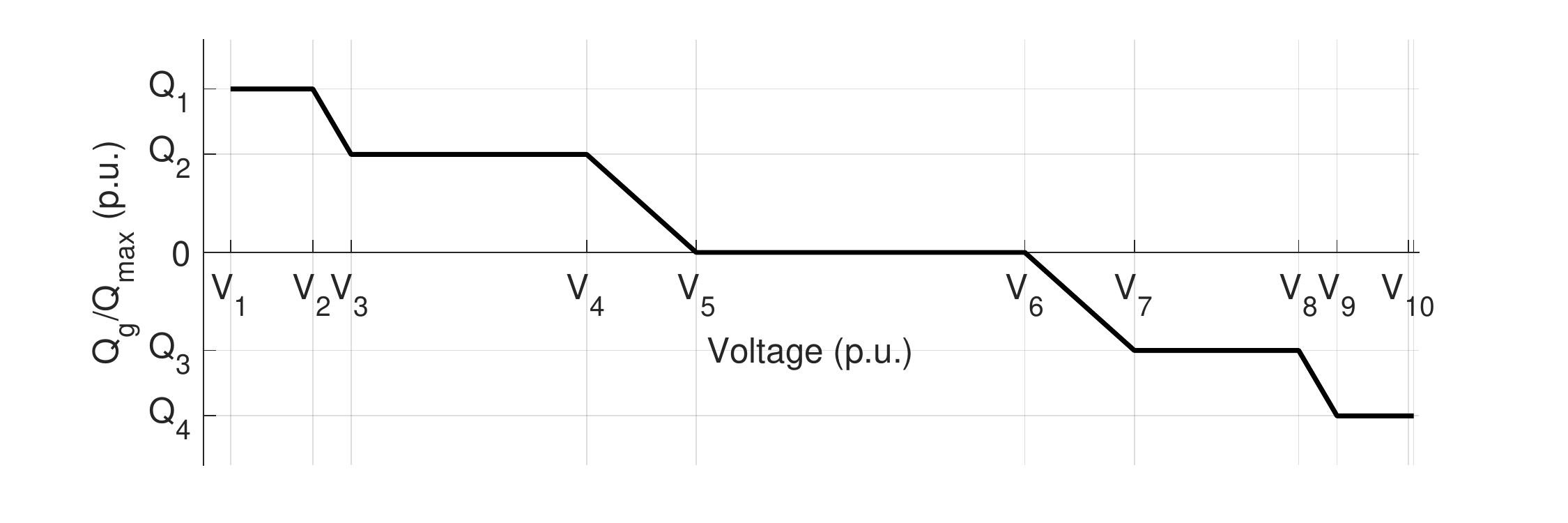}\label{char1}}
      \subfigure[]{\includegraphics[width=0.8\textwidth]{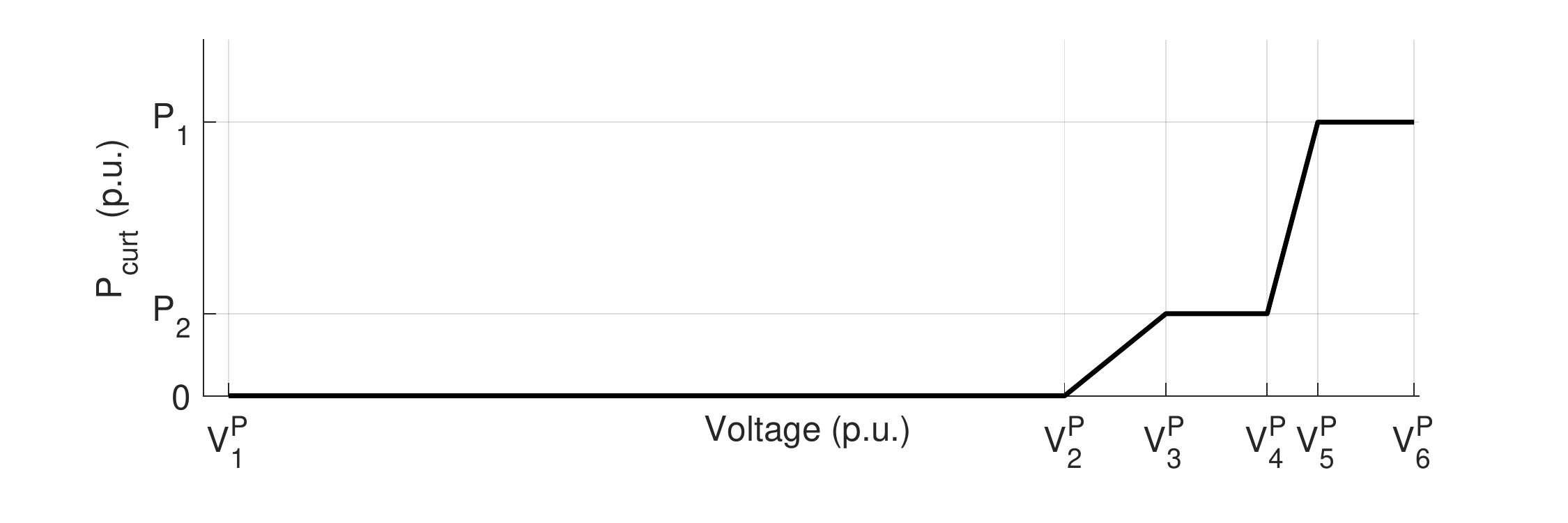}\label{char2}}
      \caption{  a) Parametric $Q-V$ characteristic.\\ b) Parametric $P_{curt}-V$ characteristic.}
       \label{char}
  \end{figure}

In this work, two characteristic curves $Q=f(V,P)$ and $P_{curt}=f(V,Q)$ are derived for each DER. First, some piecewise linear template curves are chosen for all the DERs. These curves, shown in Fig.~\ref{char}, are selected to ensure a stable response by the controller (as will be discussed later). The OPF-based solution is used to identify the curve parameters for each DER separately, while imposing some constraints that enforce the characteristics of the template curve. 
The relevant parameters and constraints are given below:
\begin{itemize}
    \item $V_1=V_1^P$ and $V_{10}=V_6^P$ are the under- and over-voltage protection limits, respectively, that when violated the DER disconnects. These are the same for all DERs.
	\item $Q_{1} \geq Q_{2} > 0$: at voltages $V \in (V_1, V_5)$, the inverter shows a capacitive behaviour to increase the depressed voltage.
    \item $Q_{4} \leq Q_{3} < 0$: at voltages $V \in (V_6, V_{10})$, the inverter consumes     reactive power and reduces the local voltage. 
    \item $V_{5} \leq V_{6}$: a deadband $V \in (V_5, V_{6})$ can be defined, in which no reactive power control is implemented.  
    \item $P_{1} \geq P_{2} > 0$: these values define the levels of curtailment when $V \in (V_2^{P}, V_{6}^{P})$ 
    \item $V_2^{P} \geq V_{9}$: to prioritize the use of reactive power before curtailing active power.
\end{itemize}

Appropriate selection of the parameters leads to the conventional curves sketched in Figs.~\ref{localschemes}-\ref{localschemes_P}. However, our approach allows for higher flexibility with multi-step curves and customized values for each DG, extracted from the OPF solutions.

\textcolor{black}{The procedure to derive the characteristic curves is summarized in Algorithm~\ref{alg_curves} and detailed here. First, vectors $S_1$ and $S_2$ are formed, comprising respectively of pairs $s_1(t)$ and $s_2(t)$ of the optimal reactive power injection and active power curtailment with the corresponding OPF voltage for each time step $t$. Then, these pairs are sorted in ascending order with respect to the voltage, and truncated to keep only the ones corresponding to time instances when the DER active power injection is higher than a predetermined threshold $P_{\textrm{thr}}$~(Step 1). Afterwards, vectors $S_{1}^{'}$ and $S_{2}^{'}$ are formed corresponding to the final characteristic curves. This is done by going through the vectors $\overline{S_1}$, $\overline{S_2}$ starting from the first pair with capacitive behaviour, i.e. with $max(Q_{\textrm{g}}(t))\geq Q_{\textrm{thr,ind}}\geq 0$ or the first sorted pairs (Step 2). Then, the succeeding sorted pairs are included, while ensuring that the curve $Q=f(V,P)$ is non-increasing and $P_{curt}=f(V,P)$ non-decreasing (Steps 3-6). Finally, some post-processing is performed to ensure the correct behaviour of the units. For example, some predefined values are added so that the form of the final characteristic curves are according to Fig.~\ref{char} (Steps 8-12), the slope is changed when the piece-wise linear segment between two successive pairs is too steep (Steps~13-16), etc}.

\def\NoNumber#1{{\def\alglinenumber##1{}\State #1}\addtocounter{ALG@line}{-1}}

\newcommand{\ALOOP}[1]{\ALC@it\algorithmicloop\ #1%
  \begin{ALC@loop}}
\newcommand{\ENDALOOP}{\end{ALC@loop}\ALC@it\algorithmicendloop}
\renewcommand{\algorithmicrequire}{\textbf{Input:}}
\renewcommand{\algorithmicensure}{\textbf{Output:}}
\newcommand{\algorithmicbreak}{\textbf{break}}
\newcommand{\BREAK}{\STATE \algorithmicbreak}
\makeatother

\algnewcommand{\IIf}[1]{\State\algorithmicif\ #1\ \algorithmicthen}
\newcommand{\EndIIf}{\unskip\ \algorithmicend\ \algorithmicif}
\newcommand*{\skipnumber}[2][1]{%
  {\renewcommand*{\alglinenumber}[1]{}\State #2}%
  \addtocounter{ALG@line}{-#1}}
  
\begin{algorithm}[t] 
\begin{algorithmic}[1]
\addtocounter{ALG@line}{-1}
\Require $S_1= [(Q_{\textrm{g}}(t)$, $V_{\textrm{OPF}}(t))]:=[s_1(1), s_1(2), ..., s_1(n)] $, \\ \quad $S_2= [(P_{\textrm{curt}}(t)$, $V_{\textrm{OPF}}(t))]:=[s_2(1), s_2(2), ..., s_2(n)]$ $ \forall t$  
\Ensure $Q=f(V,P)$ and $P_{curt}=f(V,P)$ inverter characteristic $\rightarrow$ $S_{1}^{'}$, $S_{2}^{'}$
\State Sort $S_1$, $S_2$ in ascending voltage order for each time step $t$ if $P_g(t)>P_{\textrm{thr}}$, $\rightarrow$ $\overline{S_1}$, $\overline{S_2}$
\State \textbf{Initialize:}  $S_{1}^{'}= 
\begin{cases}
    \overline{s_1}(1),& \text{if } max(Q_{\textrm{g}}(t))\leq Q_{\textrm{thr,ind}}\\
    max(s_1(:)),              & \text{otherwise}
\end{cases} , \quad \quad \quad S_{2}^{'}=\overline{s_2}(1)$
\For{$i \in [2,3,...n]$}
    \IIf{$Q_{\textrm{g}}(\overline{s_1}(i+1)) \leq Q_{\textrm{g}}(\overline{s_1}(i))$}
    	$S_{1}^{'} = S_{1}^{'} \cup \overline{s_1}(i)$
    \EndIIf
    \IIf{$P_{\textrm{curt}}(\overline{s_2}(i+1)) \geq P_{\textrm{curt}}(\overline{s_2}(i))$}
    	 $S_{2}^{'} = S_{2}^{'} \cup {\overline{s_2}(i)} $
	\EndIIf
\EndFor

\State Add missing pairs so that the characteristic curves are complete (see~Fig.~\ref{char})
\IIf{$V_{\textrm{OPF}}(s_1^{'}(1)) \geq V_3$}
     $S_{1}^{'} = S_{1}^{'} \cup {(Q_1,V_1), (Q_1,V_2), (Q_{\textrm{g}}(s_1^{'}(1)),V_3)} $
\EndIIf
\IIf{$Q_{\textrm{g}}(s_1^{'}(end)) \geq 0$}
     $S_{1}^{'} = S_{1}^{'} \cup {(0,V_{\textrm{OPF}}(s_1^{'}(end))), (0,V_6)} $
\EndIIf
\IIf{$V_{\textrm{OPF}}(s_1^{'}(end)) \leq V_8$}
     $S_{1}^{'} = S_{1}^{'} \cup {(Q_{\textrm{g}}(s_1^{'}(end)),V_8), (Q_4,V_9), (Q_4,V_{10})} $
\EndIIf
\IIf{$V_{\textrm{OPF}}(s_2^{'}(1)) \geq V_2^P$}
     $S_{2}^{'} = S_{2}^{'} \cup {(0,V_1^P), (0,V_2^P)} $
\EndIIf
\IIf{$P_{\textrm{curt}}(s_2^{'}(end)) \leq V_4^P$}
     $S_{2}^{'} = S_{2}^{'} \cup {(P_{\textrm{curt}}(s_2^{'}(end)), V_4^P), (P_1, V_5^P), (P_1, V_6^P)} $
\EndIIf
\For{$i \in [2,3,...n]$}
    
    \IIf{$slope(\overline{s_1}(i\hspace{-0.1cm} + \hspace{-0.1cm} 1),\overline{s_1}(i)) \leq tan(\phi_{\textrm{thr}})$}
    $ V_{\textrm{OPF}}(\overline{s_1}(i\hspace{-0.1cm} + \hspace{-0.1cm} 1)) -= |\frac{Q_{\textrm{g}}(\overline{s_1}(i+1)) - Q_{\textrm{g}}(\overline{s_1}(i))}{tan(\phi_{\textrm{thr}})}|$
    \EndIIf
    
    \IIf{$slope(\overline{s_2}(i\hspace{-0.1cm} + \hspace{-0.1cm} 1),\overline{s_2}(i)) \geq tan(\phi_{\textrm{thr}})$}
    $ V_{\textrm{OPF}}(\overline{s_2}(i\hspace{-0.1cm} + \hspace{-0.1cm} 1)) -= |\frac{P_{\textrm{curt}}(\overline{s_2}(i+1)) - P_{\textrm{curt}}(\overline{s_2}(i))}{tan(\phi_{\textrm{thr}})}|$
    \EndIIf
\EndFor
\\
\Return $S_{1}^{'}$, $S_{2}^{'}$
\end{algorithmic}
\caption{Derivation of $Q=f(V,P)$ and $P_{curt}=f(V,P)$ inverter characteristics} \label{alg_curves}
\end{algorithm}

The final $Q=f(V,P)$ and $P_{curt}=f(V,P)$ characteristics represent closed-loop control schemes and some stability problems might occur. For example,  if the slope of the curve is too steep or the network impedance is high, then some oscillatory behaviour might be observed. This phenomenon is analysed in~\cite{Farivar2013,Bucher2014PhD} %phdselida61-39
and~\cite{Stetz2014}%3.3selida41
, along with possible ways to mitigate the problem, such as a stabilizer or an iterative procedure with a damping term.
\textcolor{black}{In this work, we also included an iterative procedure which stops once the voltage difference between two subsequent iterations is smaller than a predetermined value.}

\textcolor{black}{The proposed characteristics can be discretized and implemented through a look-up table using a simple  microprocessor (e.g., Arduino) installed at the device and driving the primary controller. Thus, existing inverters should be retrofitted for that purpose, but the primary controllers do not need to be replaced.
The characteristic curves of the inverters can be updated frequently, if the communication infrastructure allows it. For example, in a simplistic case, we can have four characteristic curves per year and DG, i.e. one per season. By sending signals  over power lines, we can easily switch among the different characteristics in order to derive a behaviour close to the optimal for the whole year, being able to account for seasonalities in terms of the DG injections and loading.
Evaluation of yearly results will allow us not only to quantify the yearly losses incurred by our method and compare them against the ideal OPF approach, but also the operational APC and RPC costs to be compared with the conventional grid expansion measures.\\
In case of a centralized scheme, a bi-directional communication scheme would have to be put in place along with remote terminal units at the inverters (to receive commands, decrypt them, and drive the primary controller; and, to get measurements, encrypt them, and send them to centralized unit) and the central computational unit.}

\section{Case Studies}\label{sec:CaseStudies}

\subsection{Modified LV Cigre benchmark grid}

\begin{figure}
	\centering
	\includegraphics[width=0.8\textwidth]{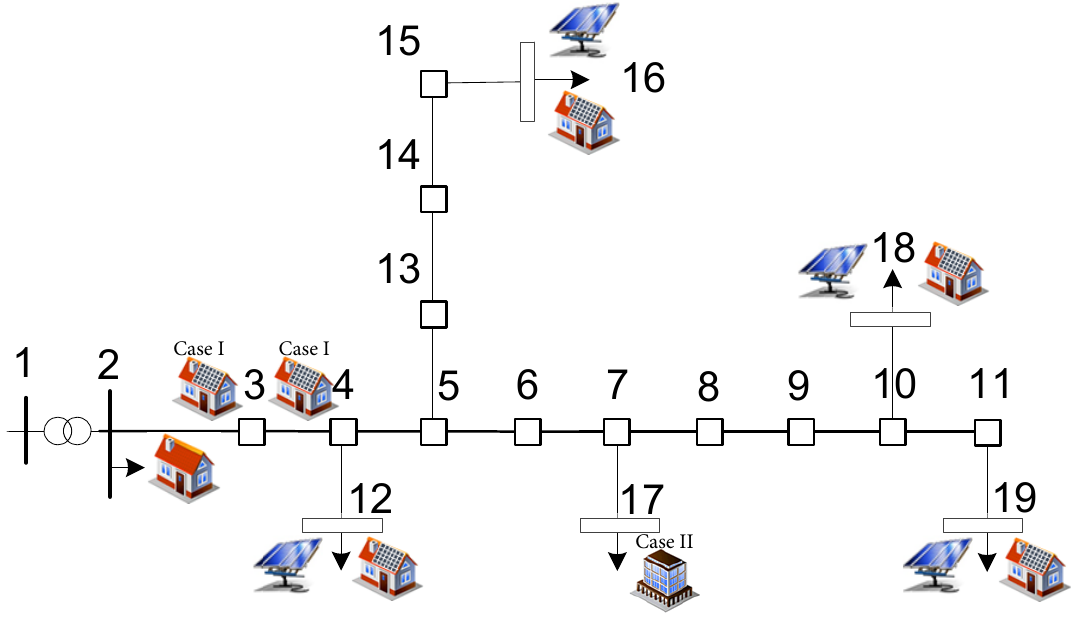} %width=1\textwidth
	\caption{LV Cigre benchmark grid - \cite{Strunz2014} modified}\label{fig:cigre_test_system}
\end{figure}

As basis for the case studies, we use the benchmark LV grid presented in~\cite{Strunz2014} and sketched in Fig.~\ref{fig:cigre_test_system}. In order to investigate future challenges of DNs, two scenarios are analysed. First, the system is modified by adding PV generators at four nodes (3, 4, 12, 16, 18, and 19). This leads to a grid facing over-voltage and transformer overload issues and their mitigation is analysed in Case~1. Then, in addition to the PV generators at nodes (12, 16, 18, and 19), the load at Node~17 is increased, representing a commercial centre. This change leads to a system with both over- and under-voltage problems, and their mitigation is studied in Case~2. In all cases, the thermal limits of the cables are increased accordingly, in order to investigate the aforementioned scenarios. The system parameters are summarized in Appendix~\ref{appendix}.

\subsubsection{Load and generation profiles}
Daily load profiles for the residential and commercial loads can be found in~\cite{Strunz2014}, in the form of scaling factors. However, the seasonal characteristics are not provided, thus the procedure described in~\cite{V.Poulios2015} is used to produce load and PV injection data. 

In order to evaluate the yearly effect of the proposed methodology, we have to consider the response of the active units under different load and DG infeed conditions. In our case, we consider one typical representative
and one worst-case, in terms of PV injection, daily profile for each season, and we assume that we can switch among 4 different characteristic curves, one for each season. Thus, we create the curves based on the 4 worst days and use them to evaluate the response of the typical days. For a yearly simulation, we assume we have 10 worst and 80 typical days per season. In case the derived curves are incomplete, in our simulations we use $P_{\textrm{thr}}=0.1 \textrm{p.u.}$, $Q_{\textrm{thr,ind}}=0.1$p.u., $\phi_{\textrm{thr}}$=88.85$^{\circ}$, $V_1=0.88\textrm{V}, V_2=0.9\textrm{V}, V_3=0.92\textrm{V}, V_5=V_{\textrm{OPF}}(s_1^{'}(end)), V_6=1.079\textrm{V}, V_8=1.08\textrm{V}, V_9=1.1\textrm{V},V_{10}=1.3\textrm{V}, Q_1=1\textrm{p.u.}, Q_2=Q_{\textrm{g}}(s_1^{'}(1)), Q_3=Q_{\textrm{g}}(s_1^{'}(end))\textrm{p.u.}, Q_4=-1\textrm{p.u.}, V_1^P=0.88\textrm{V}, V_2^P=0.92\textrm{V}, V_4^P=1.1\textrm{V}, V_5^P=1.12\textrm{V}, V_6^P=1.3\textrm{V}$ and $P_1=1\textrm{p.u.}$,  referring to steps 8-12 of Algorithm~\ref{alg_curves}.

For each case study, four different investigations are performed:
\begin{itemize}
	\item Method 0: For each time step, an AC PF solution of the LV grid based on the forecasted load consumption and PV generation is performed. The DERs are operating having a power factor of one (i.e., without reactive control);
    \item Method 1: An AC OPF solution is performed for each time step, allowing active power curtailment and reactive power control with $cos\phi_{max}=0.9$ (see Section~\ref{inverters});
    \item Method 2: An AC PF solution is performed for each time step, with the PV generators behaving according to the characteristic curve of the German grid-codes~\cite{VDE}; and 
    \item Method 3: An AC PF solution is performed for each time step, with the PV generators behaving according to the characteristic curves of the proposed hybrid approach (extracted from the results of Method~1).
\end{itemize}

\subsection{Case 1 - Over-voltage and transformer over-load problems}

The first case investigates the situation where, if no measures are taken, the transformer is overloaded and the voltage at one node exceeds its maximum allowable value. First, the worst day of the year in terms of PV production is analyzed in more detail, and then a summary of the yearly behaviour of the proposed method provides inputs for the planning decisions. 
\subsubsection{Worst summer day results}
\begin{figure}
  \centering
      \subfigure[]{\includegraphics[width=0.7\textwidth]{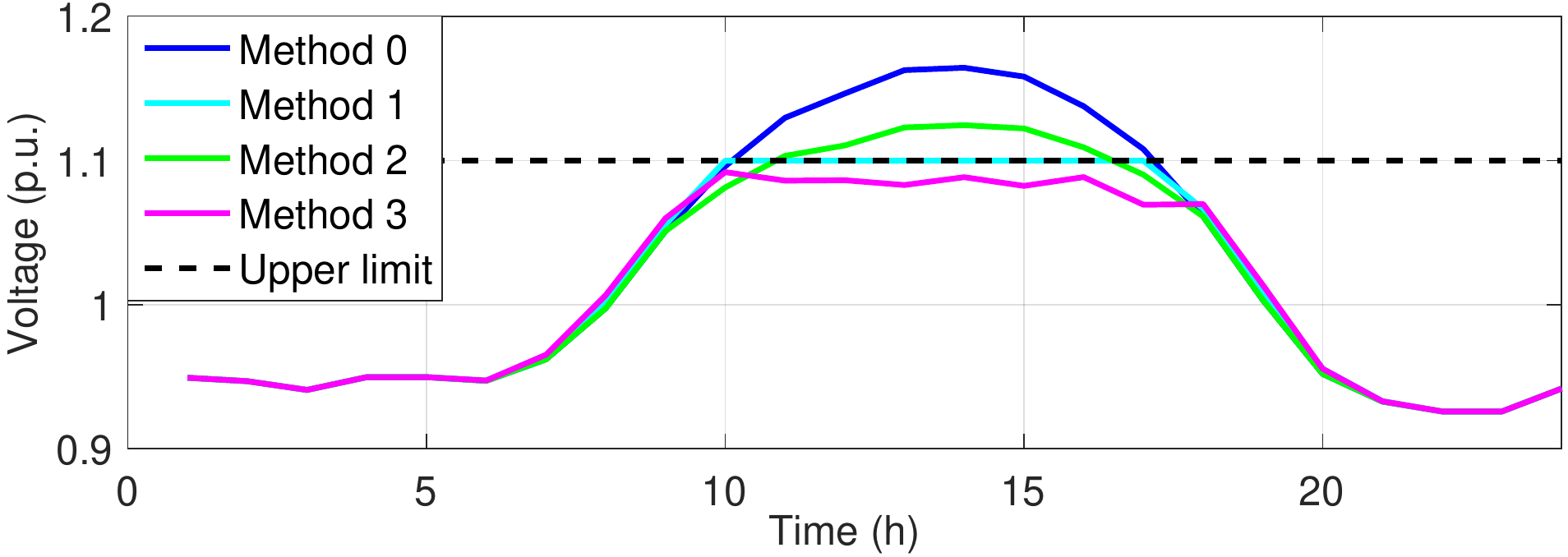}\label{Case11_V16}}
      \subfigure[]{\includegraphics[width=0.7\textwidth]{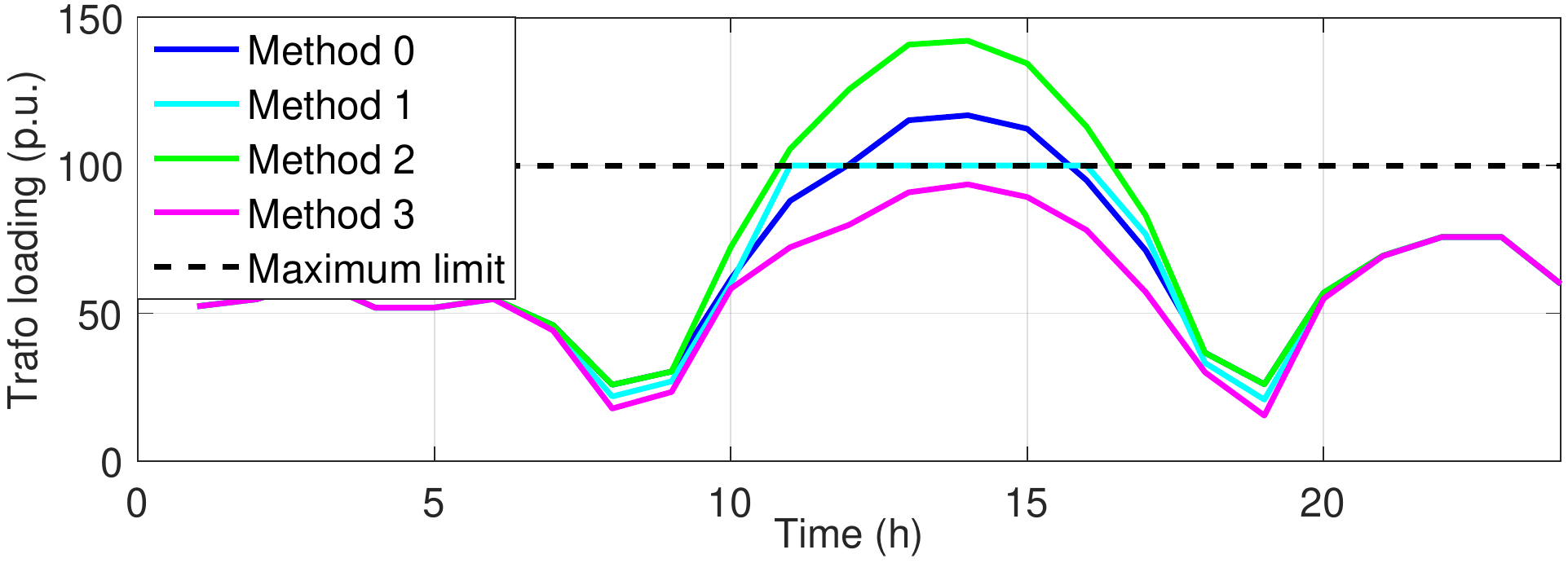}\label{Case11_TRAFO}}
      \subfigure[]{\includegraphics[width=0.7\textwidth]{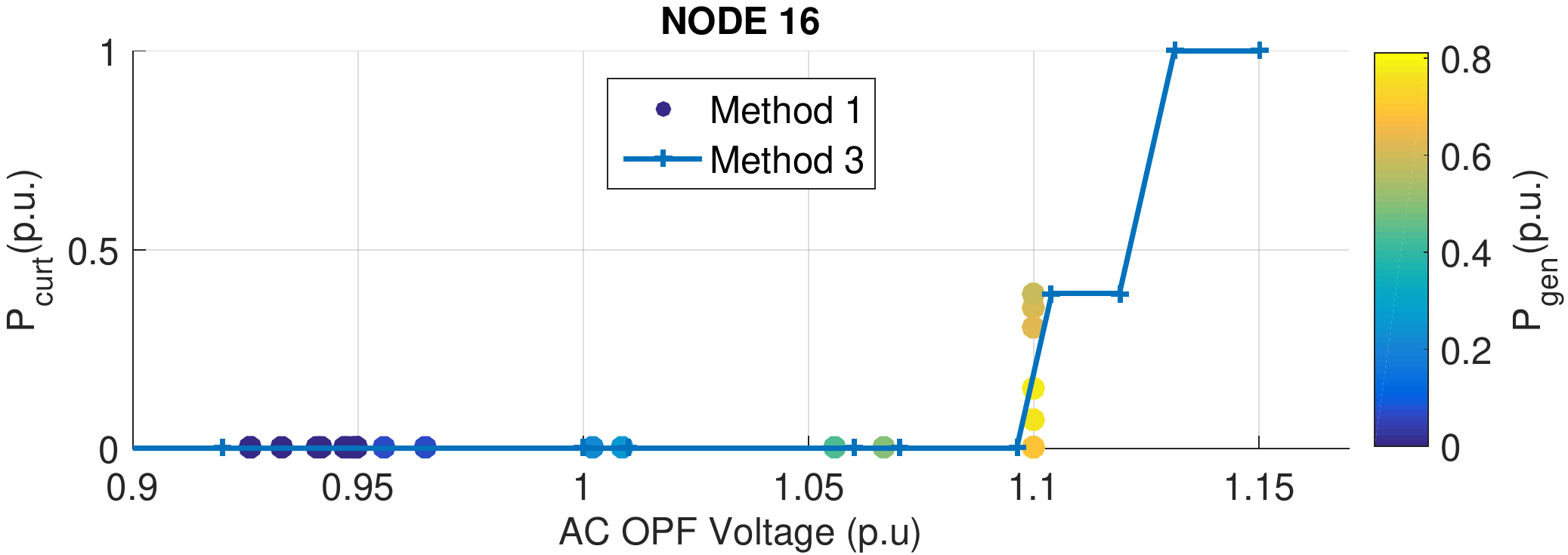}\label{Case11_Pcurt_fit}}
      \subfigure[]{\includegraphics[width=0.7\textwidth]{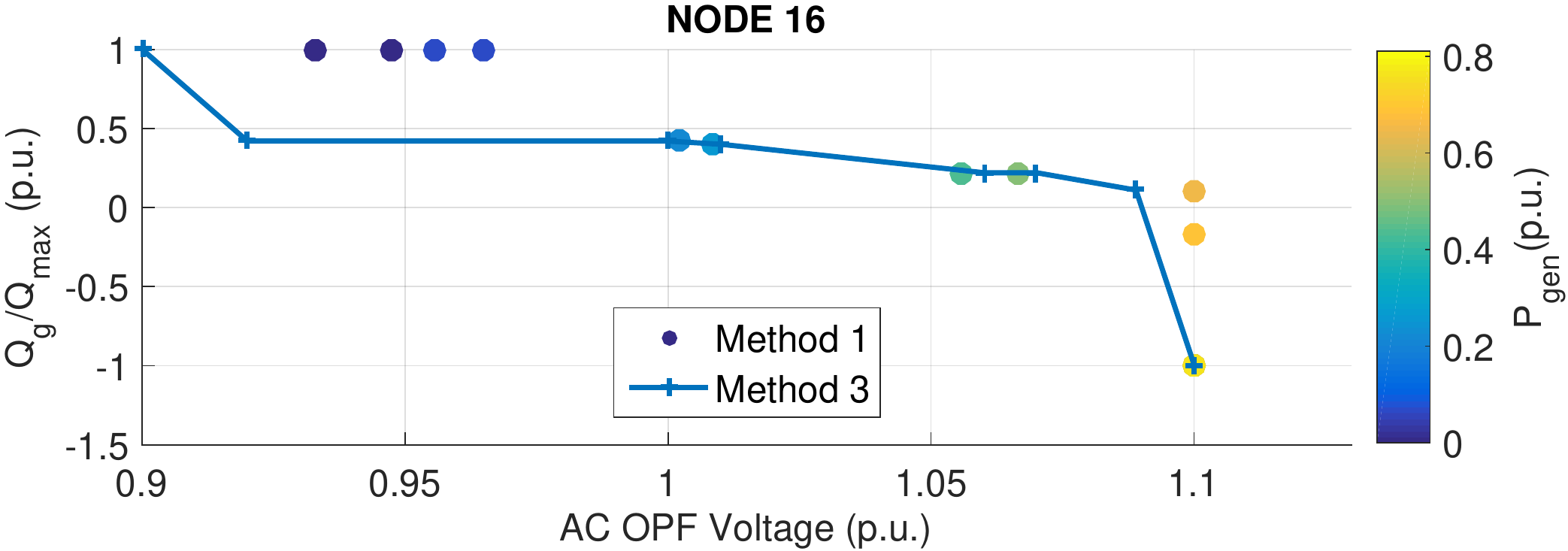}\label{Case11_QVP__fit}}
      \caption{  a) Case 1: Voltage profile at Node 16.\\ b) Case 1: Transformer loading. \\ c) Case 1: Extracted $P_{curt}=f(V)$ curve at Node 16 (Method 3).\\ d) Case 1: Extracted $Q=f(V)$ curve at all at Node 16 (Method 3).}
  \end{figure}

Figure~\ref{Case11_V16} shows the daily voltage profile
at the problematic node for all the methods considered. First, without any control (Method~$0$), the voltage rises up to 1.16 p.u. (above the maximum
acceptable value of 1.10 p.u.), indicating that measures need to be taken. Controlling the PVs
according to the German grid-code (Method~$2$) is not adequate to bring the voltage back to normal
values for the existing parameters and the assumed PV installed capacity. The OPF-based control
of the PVs (Method~$1$) satisfies the voltage constraints by using reactive power control and active
power curtailment in an optimal way. Finally, the hybrid method (Method~$3$) achieves almost the
same results as Method 1, but in a sub-optimal way as will be discussed below.

The loading of the transformer is compared for all methods in Fig.~\ref{Case11_TRAFO}. Method~$2$ shows the maximum loading, even higher than Method~$0$, due to the increased reactive power needs from the PV nodes. Method~$1$ satisfies the constraint at its limit, whereas Method~$3$ does not lead to thermal violations.   
 
Figure~\ref{Case11_Pcurt_fit} shows the $P_{curt}=f(V,P)$ characteristic curve of Node $16$ for Method~3, as extracted from the OPF-based solution. The goal here is to alleviate not only the over-voltage, but also the overload of the transformer. Figure~\ref{Case11_QVP__fit} shows the $Q=f(V,P)$ characteristic curve at Node~$16$, which is interesting since in the optimal response it behaves both capacitive and inductive at different hours. At low voltages, this node is injecting reactive power to reduce the losses of the system, while at voltages higher than $1.085$ p.u., it becomes inductive to reduce the local voltage. This is not captured by the current characteristic curves which consider only the production of reactive power, occurring at high active power injection.
 
\begin{figure}
 	\centering
 	\subfigure[]{\includegraphics[width=1\textwidth]{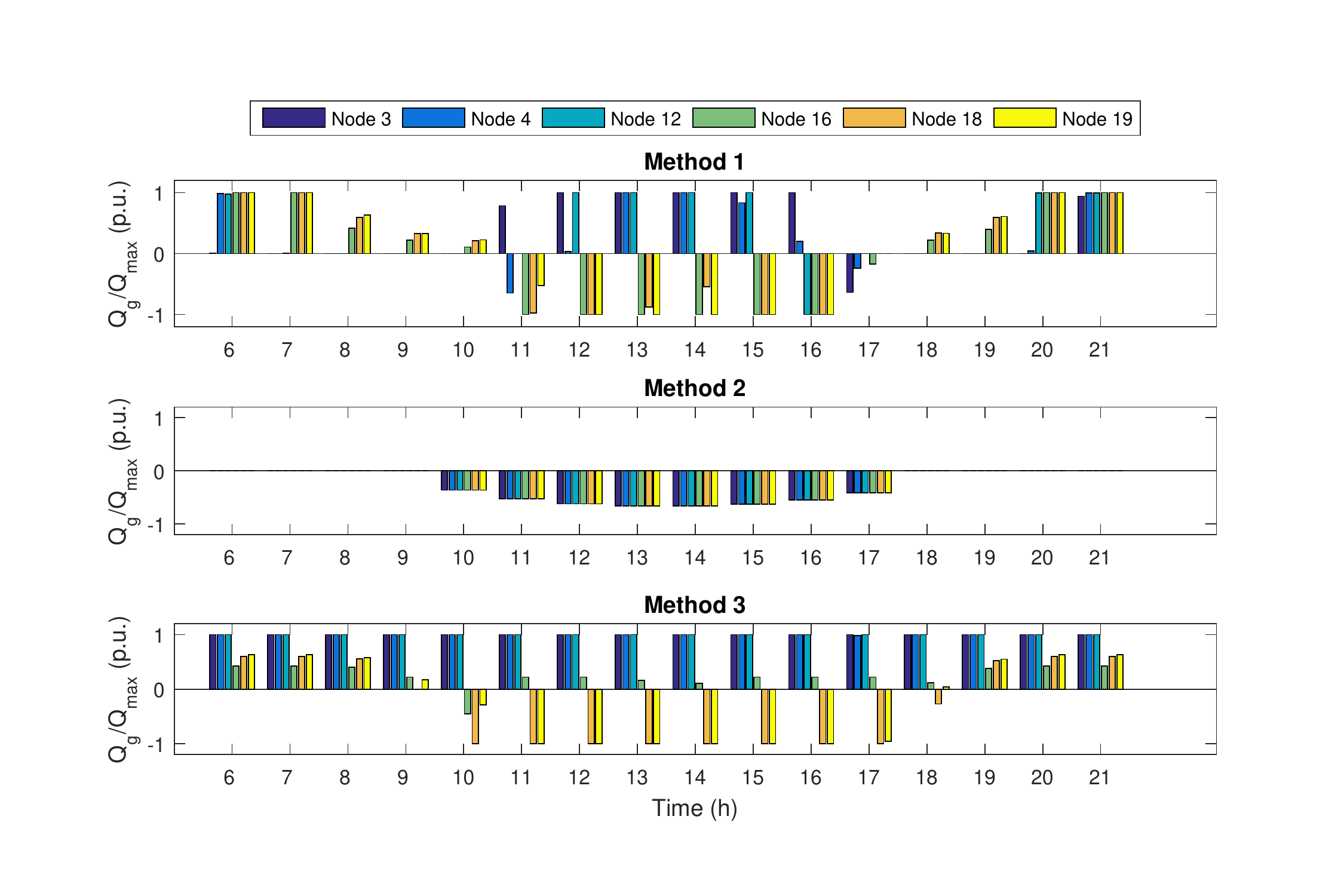}\label{Case11_Qcon}}
    \subfigure[]{\includegraphics[width=0.7\textwidth]{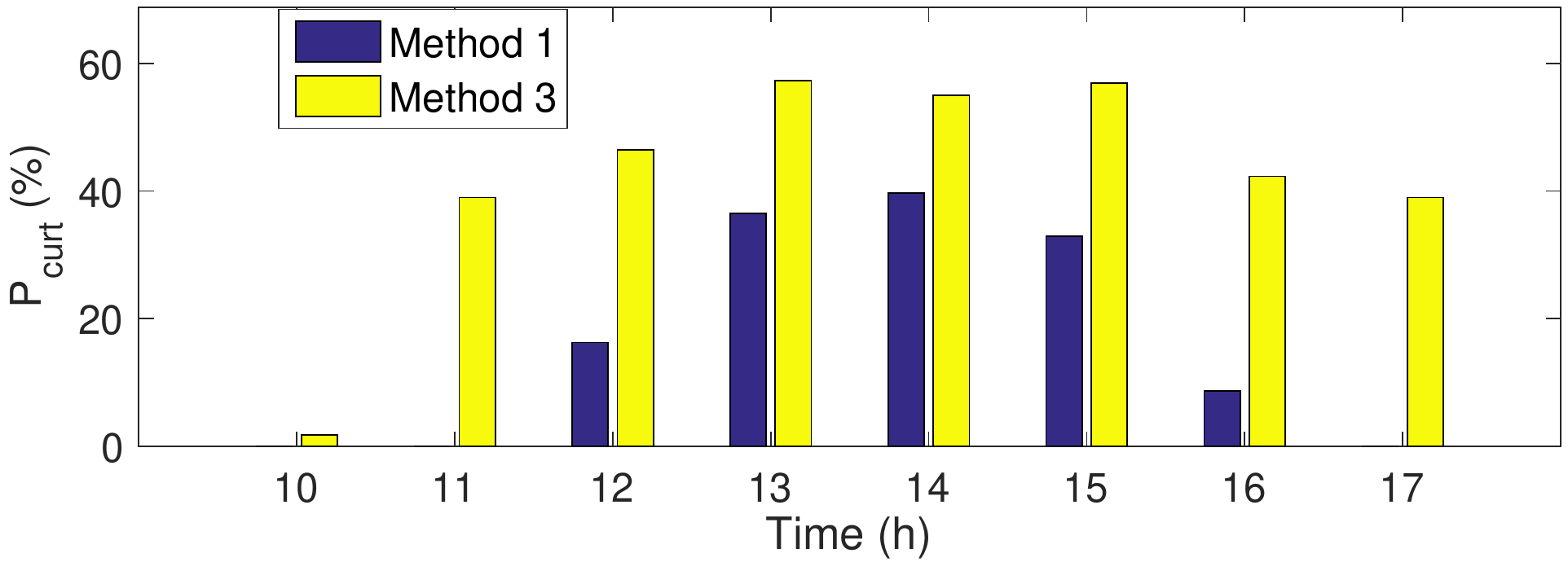}\label{Case1_Pcurt}}
     	\caption{a) Case 1: Reactive power control for all configurations in relevant hours (zero otherwise) \\ b) Case 1: Curtailed active power at Node 16 in relevant hours (zero otherwise)
        }	\label{Case11_Qcon}
 \end{figure}

Figure~\ref{Case11_Qcon} shows the use of reactive power control for the different methods for the relevant hours. According to Method $1$, Nodes~$3$, $4$ and $12$, show a capacitive behaviour at problematic hours~12:00-15:00. In this way, they produce a part of the needed reactive power for the whole grid, reducing the requirement of the transformer. The other nodes consume reactive power, e.g. Node $16$, to reduce its local voltage. Method~$2$ results in consumption of reactive power at all PV nodes by definition. Since a perfect spacial correlation of the PV infeed is assumed, all nodes are consuming the same p.u. reactive power. Thus, on one hand voltage drops further in all nodes, but on the other hand, the transformer should carry the increased reactive power needs. Method~$3$ tries to emulate the behaviour of Method~$1$, as seen also in Fig.~\ref{Case1_Pcurt} which compares the power curtailment of Methods~$1$ and ~$3$. Method~$3$ results in increased curtailment, even at hours when curtailment is not needed to satisfy the constraints.

Finally, Table~\ref{Case11res} summarizes the results of this day. It can be seen that all Methods satisfy the lower voltage limit, but only Methods~$1$ and~$3$ the upper limit. Methods~$1$ and~$3$ require PV curtailment with the latter exhibiting higher active power curtailment. Therefore, a comparison of the losses is not meaningful here, but overall Method~$2$ shows the highest, due to increased reactive power needs. This results also in the largest transformer over-loading. Method~3 does not violate either of the voltage or transformer limits.  

\begin{table}
\centering
\caption{Case 1: Summarized worst-summer day results for all methods.}
\label{Case11res}
\begin{tabular}{c|cccc}
                          & Method 0 & Method 1 & Method 2 & Method 3 \\ \hline
$V_{max}$ (p.u.)             & $\textcolor{red}{1.1645}$   & $1.1000$   & $\textcolor{red}{1.1246}$   & $\textcolor{orange}{1.0921}$   \\
$V_{min}$ (p.u.)             & $0.9259$   & $0.9259$   & $0.9259$   & $0.9259$   \\
Losses (\%)               & $7.1027$   & $7.1588$   & $8.7069$   & $5.4578$   \\
$P_{curt} (\%)$           & $0$        & $3.87$   & $0$        & $8.80$   \\
Maximum transformer loading (\%)  & $\textcolor{red}{117.00}$    & $100.00$      & $\textcolor{red}{142.22}$    & $\textcolor{orange}{93.63}$      \\
Security constraints satisfied & \textcolor{red}{NO}       & \textcolor{green}{YES}      & \textcolor{red}{NO}       & \textcolor{orange}{YES}       \\ \hline
\end{tabular}
\end{table}

\subsubsection{Yearly evaluation - Planning recommendations}
The yearly evaluation is summarized in Table~\ref{yearly11}, where the results of Method~1 are used as a benchmark. The results for maximum voltage occur during the worst spring day and are similar to the worst summer day. There are no under-voltage issues, since the lowest voltage occurs at the typical winter day. Furthermore, we observe that Method~2 shows significantly higher losses compared to the OPF solution. The proposed method shows overall reduced losses, i.e. $4.85\%$ more than Method~1. However, this is explained by the fact that it leads also to higher power curtailment. Thus, the losses comparison should be always linked with curtailment results as well as the constraint satisfaction capability, so that the conclusions are not misleading. Method~2 shows the worst maximum transformer loading due to additional reactive power needs, while Method~3 shows no overload. 

The conventional planning procedure would result in ways to increase the capacity of the transformer for Methods~$0$ and $2$, having an impact on the capital expenditure of the DSO. By using Method $3$, co-optimizing the stages of planning and operation, the DSO can use the existing infrastructure, deferring investments, with marginal violations, which might be acceptable for limited time, i.e. a few hours every year.  

\begin{table}
\centering
\caption{Case 1: Summarized yearly results for all methods.}
\label{yearly11}
\begin{tabular}{c|cccc}
                          & Method 0 & Method 1 & Method 2 & Method 3 \\ \hline
$V_{max}$ (p.u.)             & \textcolor{red}{$1.1649$}    & $1.1000$   & \textcolor{red}{$1.1246$}   & \textcolor{orange}{$1.1017$}   \\
$V_{min}$ (p.u.)            & $0.9029$   & $0.9029$   & $0.9029$   & $0.9029$   \\
Relative losses with respect to Method 1 $(\%)$               & $+ 1.6$   & $1$ (ref)  & $+ 7.28$   & $- 4.85$   \\
Relative $P_{curt}$ with respect to Method 1 $(\%)$           & --        & $1$ (ref)  & --        & $+4.73$   \\
Maximum transformer loading $(\%)$  & \textcolor{red}{$117.00$}    & $100.00$      & \textcolor{red}{$144.50$ }     & $98.50$      \\
Security constraints satisfied & \textcolor{red}{NO}       & \textcolor{green}{YES}      & \textcolor{red}{NO}       & \textcolor{orange}{YES}       \\ \hline
\end{tabular}
\end{table}

\subsection{Case 2 - Over- and under-voltage problems}
In this case, the grid is facing simultaneously an over-voltage issue at Node $16$, and under-voltage problems at the commercial load terminals (Node~$17$). It should be noted that there is no control available at Node~$17$ and thus, the suitability of the examined methods to tackle system-wide issues is evaluated.

\subsubsection{Worst summer day results}
\begin{figure}
  \centering
      \subfigure[]{\includegraphics[width=0.7\textwidth]{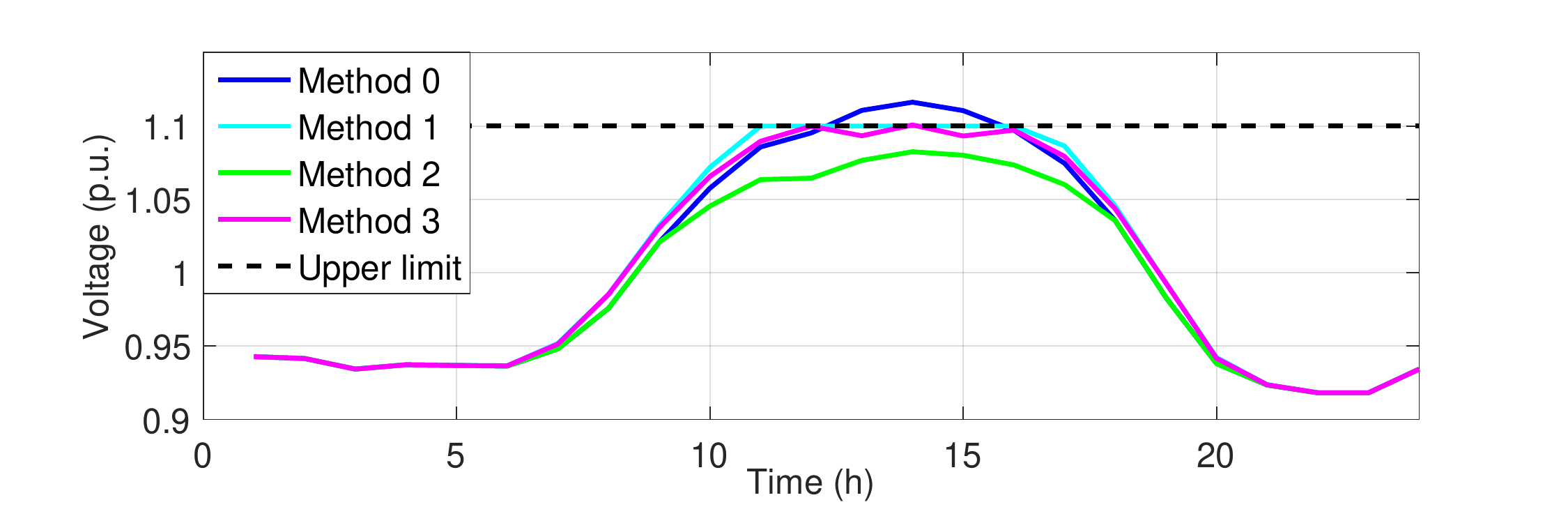}\label{Case2_V16}}
      \subfigure[]{\includegraphics[width=0.7\textwidth]{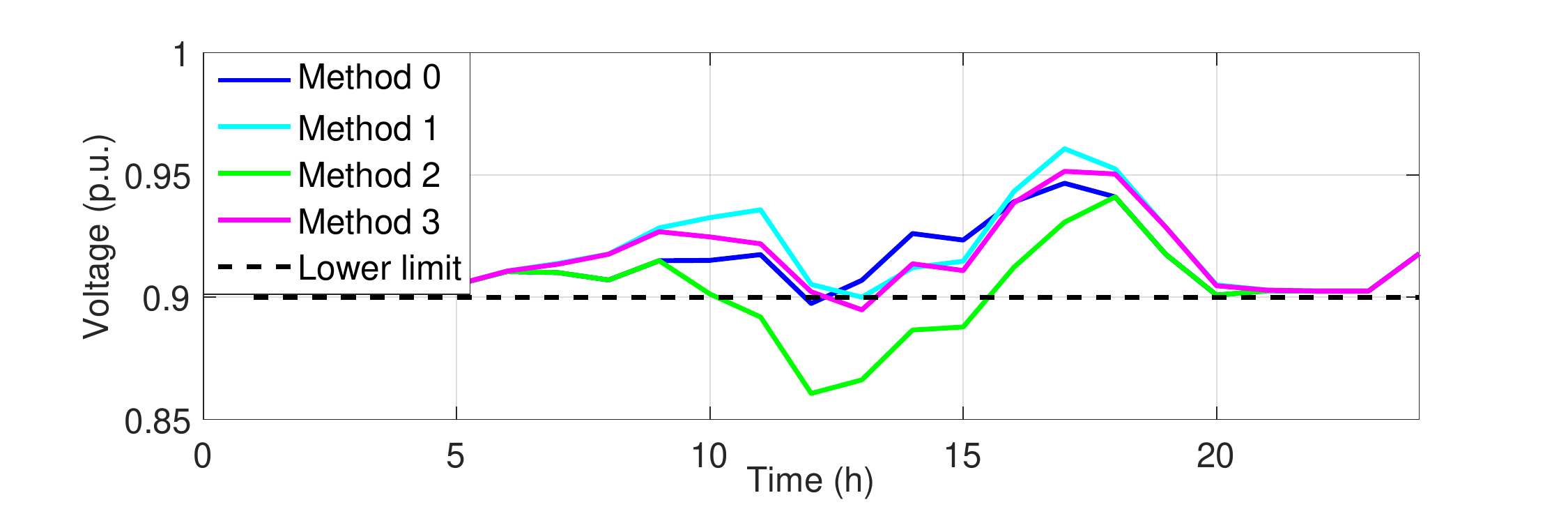}\label{Case2_V17}}
      \subfigure[]{\includegraphics[width=0.7\textwidth]{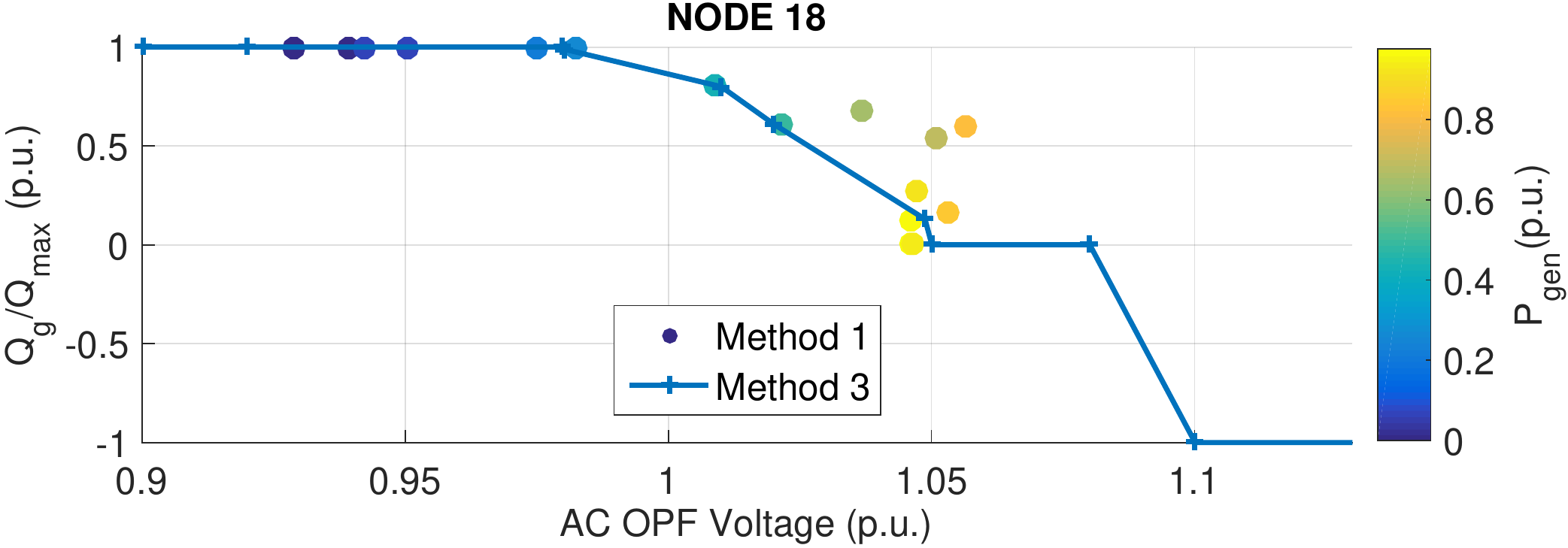}\label{Case2_Pcurt}}
      \subfigure[]{\includegraphics[width=0.7\textwidth]{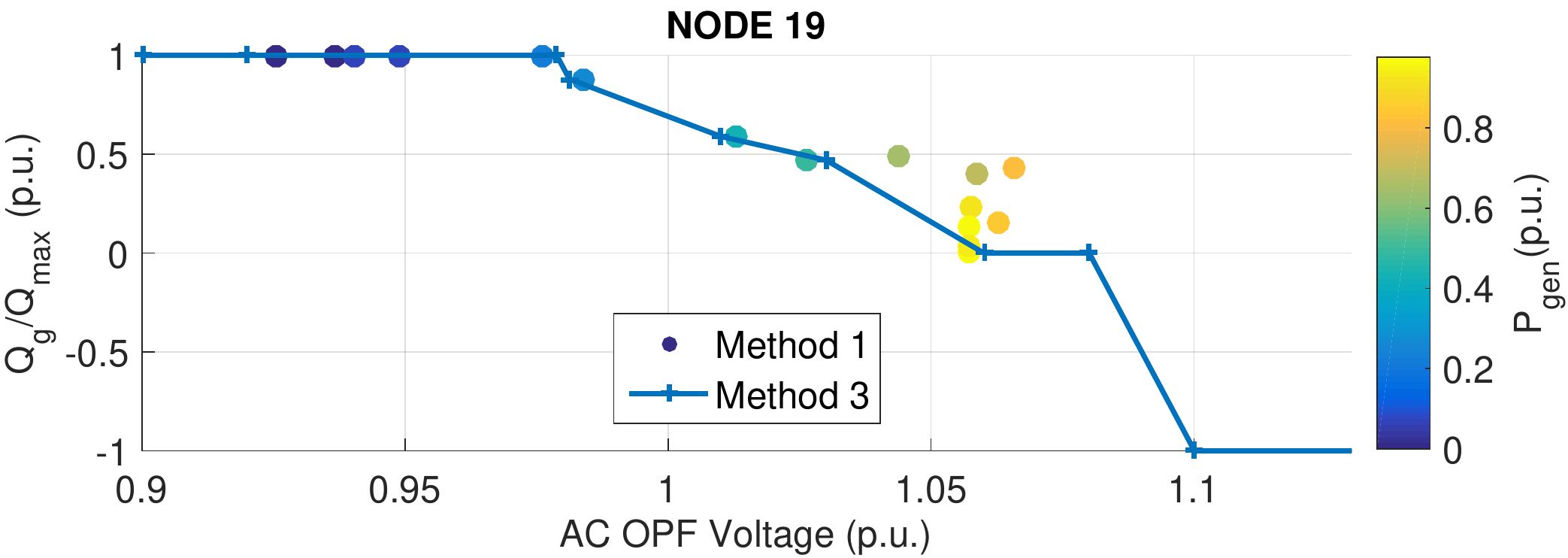}\label{Case2_QVP__fit}}
      \caption{  a) Case 2: Voltage profile at Node 16.\\ b) Case 2: Voltage profile at Node 17. \\ c)  Case 2: Extracted $Q=f(V)$ curve at PV Node 18 (Method 3).\\ d) Case 2: Extracted $Q=f(V)$ curve at PV Node 19 (Method 3).}
  \end{figure}

Figures~\ref{Case2_V16} and~\ref{Case2_V17} show the voltage profiles at both problematic nodes. As can be seen, Method~1 tackles both issues by choosing optimal set-points for the DGs. Method~2 satisfies the upper voltage constraint, but  reduces even further the voltage at the node of the commercial load, having a worsening effect in this case. Finally, Method~3 succeeds to secure the system, since marginal violations are acceptable for limited time. 

The characteristic curves of the PV at Node~$16$ are similar to Case~$1$. However, the curves of the PVs at Nodes~$18$ and~$19$ -shown in Fig.~\ref{Case2_QVP__fit}- are more interesting: a capacitive behaviour is \textcolor{black}{observed over a large range of local voltages}, to support the remote Node~$17$.

\begin{figure}
 	\centering
    {\includegraphics[width=1\textwidth]{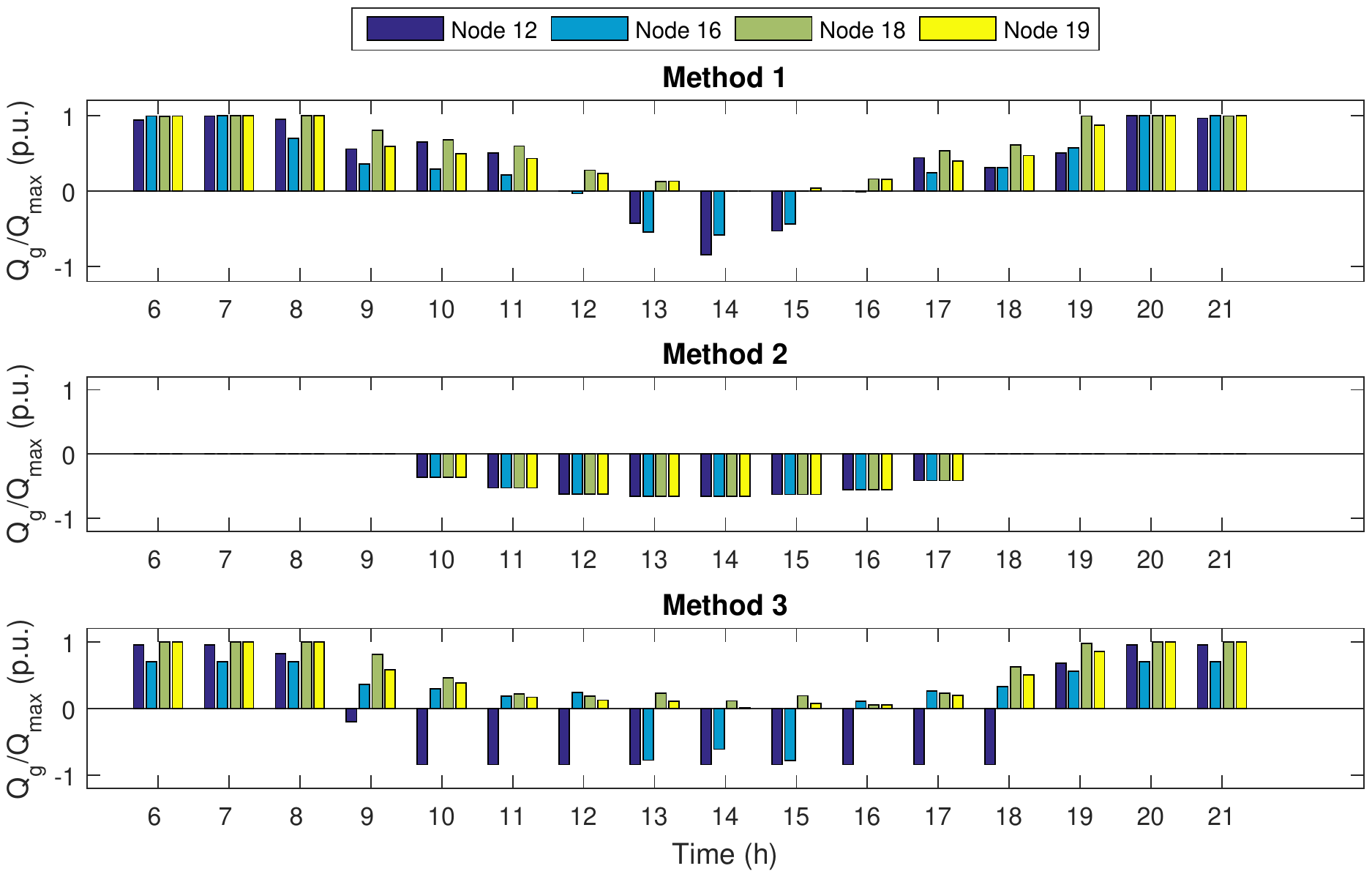}}
   	\caption{Case 2: Reactive power control for all configurations in relevant hours (zero otherwise)}	\label{Case22_Qcon}
 \end{figure}

The reactive power injections of the different methods are shown in Fig.~\ref{Case22_Qcon}. Method $1$ defines the optimal behaviour of the inverters as follows: Node~$12$ becomes clearly inductive only during three hours to assist Node~$16$ to reduce its local voltage. At the same time, Nodes~$18$ and~$19$ become capacitive to increase the voltage of Node~$17$. Method~$2$ increases the consumption of reactive power at all PV nodes by definition, and hence, voltage drops further in all nodes. On the other hand, Method~$3$ which tries to mimic the behaviour of Method~$1$, leads to marginal  voltage violations.

Finally, Table~\ref{Case22res} summarizes the results of interest. While Method~2 satisfies the upper voltage limit, it deviates from the acceptable minimum value of $0.9$ p.u. by $4.37$\%. This is expected as this control is designed only to cope with over-voltage problems and in a uniform way throughout the grid. Method~$1$ satisfies both the upper and lower voltage constraints, while Method~$3$ with marginal violations. Compared to the OPF solution, Methods~$2$ and~$3$ show increased losses, namely by~$24.67\%$ and~$2.60\%$. In this case, there is no need for active power curtailment.

\begin{table}
\centering
\caption{Case 2:  Summarized worst summer day results for all methods.}
\label{Case22res}
\begin{tabular}{c|cccc}
                          & Method 0 & Method 1 & Method 2 & Method 3 \\ \hline
$V_{max}$ (p.u.)             & \textcolor{red}{$1.1162$}   & $1.1000$   & $1.0825$   & \textcolor{orange}{$1.1007$}   \\
$V_{min}$ (p.u.)             & $0.8975$   & $0.9000$   & \textcolor{red}{$0.8607$}   & \textcolor{orange}{$0.8948$}   \\
Losses $(\%) $              & $7.4754$   & $7.3995$   & $9.2250$   & $7.5919$  \\
Total $P_{curt} (\%)$           & $0$        & $0$  	    & $0$        & $0$      \\
Security constraints satisfied & \textcolor{red}{NO}       & \textcolor{green}{YES}      & \textcolor{red}{NO}       & \textcolor{orange}{YES}      \\ \hline
\end{tabular}
\end{table}

\subsubsection{Yearly evaluation - Planning recommendations}
The yearly evaluation is summarized in Table~\ref{yearly2}. There are under-voltage issues in all Methods apart from Method~1. However, Method~3 shows only marginal violation compared with Method~2 which deteriorates the under-voltage issue. 

Since there is no active power curtailment for any of the Methods, the comparison of the losses is more meaningful in this case. As we can observe, Methods~2 and~3 result in more losses compared to the OPF solution, by namely roughly 20$\%$ and 8$\%$. The difference comes from the fact that Method~3 can trigger a capacitive behaviour at some nodes, which inject part of the needed reactive power for the loads.

The conventional planning procedure would highlight the need to take local measures,  e.g. installing shunt capacitors at Node~$17$, against under- and over-voltages. 
By using the proposed Method however, the DSO can postpone investments, if the marginal violations are acceptable according to current grid codes.

\begin{table}
\centering
\caption{Case 2: Summarized yearly results for all methods.}
\label{yearly2}
\begin{tabular}{c|cccc}
                          & Method 0 & Method 1 & Method 2 & Method 3 \\ \hline
$V_{max}$ (p.u.)            & \textcolor{red}{$1.1162$}    & $1.1000$   & \textcolor{green}{$1.0825$}   & \textcolor{orange}{$1.1106$}   \\
$V_{min}$ (p.u.)            & $0.8848$   & $0.9000$   & $0.8607$   & \textcolor{orange}{$0.8948$}   \\
Relative losses with respect to Method 1 $(\%)$               & $+ 13.99$   & $1$ (ref)  & $+ 19.95$   & $+ 8.23$   \\
Security constraints satisfied & \textcolor{red}{NO}       & \textcolor{green}{YES}      & \textcolor{red}{NO}       & \textcolor{orange}{YES}       \\ \hline
\end{tabular}
\end{table}

\section{Concluding Remarks}\label{sec:Conclusions}

In the past, network reinforcement was the main tool for tackling overvoltage, undervoltage, and congestion problems in distribution grids. However, the proliferation of DGs allows for alternative solutions exploiting the control capabilities provided by these units. In this paper, a novel methodology, bridging the operational aspects of distribution grids with the planning stage is proposed, leading to a more efficient use of the existing infrastructure and of the control capabilities of DGs. 

Based on a series of off-line forecast-based OPF solutions, the local control schemes are tuned according to the location of the DGs, the system-wide challenges of the grid, and the predicted operating conditions. Then, in real-time operation, the DGs react only to local measurements adjusting their active and reactive injections based on the resulting control curves, and without the need to communicate with each other. The performance of these optimized local controllers is comparable to the OPF-based solution and can cope with more complex, system-wide problems.

Future research will explore the addition of other active units, such as controllable loads or battery energy storage systems, to avoid curtailing active power and to increase the hosting capacity of RES in distribution grids. This will add inter-temporal constraints and make the off-line optimization more complex and computationally demanding. Finally, the regulatory and privacy topics that may arise with regard to ownership and operation of active units by DSOs, will be investigated.

\bibliographystyle{ieeetr}
\bibliography{main}

\newpage
\appendix

\section{Grid parameters}
\label{appendix}
\begin{table}[h!]
	\centering
	\caption{Grid parameters for the two case studies.}
	\label{ParTab}
	\tabcolsep=0.1cm
	\begin{tabular}{|p{0.2cm} p{1.2cm} c|ccccc|}
		\hline
		& \multirow{2}{*}{\begin{tabular}[c]{@{}c@{}}Line\\ segment\end{tabular}} & \multirow{2}{*}{meters} & \multicolumn{2}{c|}{\multirow{2}{*}{Conductor ID}} & \multicolumn{3}{c|}{Phase impedance matrix ($\Omega$/km)}                                                       \\
                                           &                                      &                         & \multicolumn{2}{c|}{}                              & A                     & B                                                  & C                                                \\  \hline
 \multicolumn{1}{|l|} {\multirow{9}{*}{\rotatebox[origin=c]{90}{UG1}}} & 2-3  							                               & $35$                      & \multirow{1}{*}{UG1 / 3-ph}      & A               & $0.287+j0.167$          & $0.121+j0.110$                                       & $0.125+j0.070$                                     \\
\multicolumn{1}{|l|}{}&3-4                                                                         & $35$                      &   NA2XY                                & B               & $0.121+j0.110$          & $0.279+j0.203$                                       & $0.121+j0.110$                                     \\
\multicolumn{1}{|l|}{}&4-5                                                                         & $35$                      &    $240$ mm$^{2}$                              & C               & $0.125+j0.070$          & $0.121+j0.110$                                       & $0.287+j0.167$                                     \\ \cdashline{4-8}
\multicolumn{1}{|l|}{}&5-6                                                                         & $35$                      & \multirow{1}{*}{UG3 / 3-ph}      & A               & $1.152+j0.458$          & $0.321+j0.390$                                       & $0.330+j0.359$                                     \\
\multicolumn{1}{|l|}{}&6-7                                                                         & $35$                      &  NA2XY                                & B               & $0.321+j0.390$          & $1.134+j0.477$                                       & $0.321+j0.390$                                     \\
\multicolumn{1}{|l|}{}&7-8                                                                         & $35$                      &   $50$ mm$^{2}$                     & C               & $0.330+j0.359$          & $0.321+j0.390$                                       & $1.152+j0.458$                                     \\ \cline{4-8}
\multicolumn{1}{|l|}{}&8-9                                                                         & $35$                      & \multicolumn{5}{c|}{Transformer parameters}                                                                                                                                        \\ \cline{4-8} 
\multicolumn{1}{|l|}{}&9-10                                                                         & $35$                      & Connection                       & $V_1$ (kV)         & $V_2$ (kV)               & $Z_{tr}$ ($\Omega$)                                           & $S_{rated}$ (kVA)                                     \\
\multicolumn{1}{|l|}{}&10-11                                                                         & $35$                      & 3-ph Dyn1                        & $20$              & $0.4$                   & $0.0032+j0.0128$                                     & $500$                                             \\ \cline{1-1} \cline{4-8}
\multicolumn{1}{|l|}{\multirow{10}{*}{\rotatebox[origin=c]{90}{UG3}}} & 4-12                                                                        & $30$                      & \multicolumn{5}{c|}{Load and PV panel parameters}                                                                                                                                  \\ 
\multicolumn{1}{|l|}{}&5-13                                                                        & $35$                      & Node                             & Load (kVA)      & Power factor          & \multicolumn{2}{c|}{\begin{tabular}[c]{@{}c@{}}PV installed capacity \\ (\% / total load)\end{tabular}} \\ \cline{4-8}

\multicolumn{1}{|l|}{}&                                                                        &                      & $3$                                & $-$             & $-$                  & \multicolumn{2}{c|}{$50/-$}                                                                                 \\

\multicolumn{1}{|l|}{}&                                                                        &                      & $4$                                & $-$             & $-$                  & \multicolumn{2}{c|}{$50/-$}                                                                                 \\

\multicolumn{1}{|l|}{}&{\multirow{-4}{*}{\rotatebox[origin=c]{0}{13-14}}}                                                                        & {\multirow{-4}{*}{\rotatebox[origin=c]{0}{$35 $  }}}                   & $2$                                & $200$             & $0.95$                  & \multicolumn{2}{c|}{$-$}                                                                                 \\
\multicolumn{1}{|l|}{}&{\multirow{-2}{*}{\rotatebox[origin=c]{0}{14-15 }}}                                                                       & {\multirow{-2}{*}{\rotatebox[origin=c]{0}{$35 $  }}}                         & $12$                               & $15$              & $0.95$                  & \multicolumn{2}{c|}{$15/15$}                                                                            \\
\multicolumn{1}{|l|}{}&15-16                                                                       & $30$                      & $16$                               & $52$              & $0.95$                  & \multicolumn{2}{c|}{$60/42$}                                                                               \\
\multicolumn{1}{|l|}{}&7-17                                                                       & $30/90$                      & $17$                               & $55/210$       & $0.95/0.85$                  & \multicolumn{2}{c|}{$-$}                                                                                 \\
\multicolumn{1}{|l|}{}&10-18                                                                        & $30$                      & $18$                               & $35$              & $0.95$                  & \multicolumn{2}{c|}{$30/30$}                                                                            \\
\multicolumn{1}{|l|}{}&11-19                                                                        & $30$                      & $19$                               & $47$              & $0.95$                  & \multicolumn{2}{c|}{$40/40$}                                                                           \\
\hline
\end{tabular}
\end{table} 

\end{document}